\documentclass{amsart}\usepackage{amsmath}

\usepackage{amssymb}
\usepackage{amscd}
\usepackage{thmdefs}

\input{tcilatex}
\theoremstyle{definition}
\theoremstyle{remark}
\numberwithin{equation}{section}

\input tcilatex

\begin{document}
\title[Compatibility And Its Application]{Compatibility of a noncommutative
probability space and a noncommutative probability space with amalgamation
And Its Application}
\author{Ilwoo Cho}
\address{Dep. of Math, Univ. of Iowa, Iowa City, IA, U. S. A}
\email{icho@math.uiowa.edu}
\maketitle

\begin{abstract}
In this paper, we will consider R-transform theory and R-transform calculus
for compatible noncommutative probability space and amagamated
noncommutative probability space (See [18]). By doing this, we can realize
the relation between scalar-valued R-transforms and operator-valued moment
series, under the compatibility. 
\end{abstract}

\strut

\strut

Voiculescu developed Free Probability Theory. Here, the classical concept of
Independence in Probability theory is replaced by a noncommutative analogue
called Freeness (See [20]). There are two approaches to study Free
Probability Theory. One of them is the original analytic approach of
Voiculescu and the other one is the combinatorial approach of Speicher and
Nica (See [19], [1] and [17]).\medskip

To observe the free additive convolution and free multiplicative convolution
of two distributions induced by free random variables in a noncommutative
probability space (over $B=\mathbb{C}$), Voiculescu defined the R-transform
and the S-transform, respectively. These show that to study distributions is
to study certain ($B$-)formal series for arbitrary noncommutative
indeterminants.\strut

\strut \strut

Speicher defined the free cumulants which are the main objects in
Combinatorial approach of Free Probability Theory. And he developed free
probability theory by using Combinatorics and Lattice theory on collections
of noncrossing partitions (See [17]). Also, Speicher considered the
operator-valued free probability theory, which is also defined and observed
analytically by Voiculescu, when $\mathbb{C}$ is replaced to an arbitrary
algebra $B$ (See [19]). Nica defined R-transforms of several random
variables (See [1]). He defined these R-transforms as multivariable formal
series in noncommutative several indeterminants. To observe the R-transform,
the M\"{o}bius Inversion under the embedding of lattices plays a key role
(See [19],[17],[5],[20],[20], [17], [31] and [32]).\strut

\strut

In [20], [19] and [20], we observed the amalgamated R-transform calculus.
Actually, amalgamated R-transforms are defined originally by Voiculescu (See
[11]) and are characterized combinatorially by Speicher (See [19]). In [20],
we defined amalgamated R-transforms slightly differently from those defined
in [19] and [11]. We defined them as $B$-formal series and tried to
characterize, like in [1] and [17]. The main tool which is considered, for
studying amalgamated R-transforms is amalgamated boxed convolution. After
defining boxed convolution over an arbitrary algebra $B,$ we could get that

\strut

\begin{center}
$R_{x_{1},...,x_{s}}\,\,\frame{*}\,%
\,R_{y_{1},...,y_{s}}^{symm(1_{B})}=R_{x_{1}y_{1},...,x_{s}y_{s}},$ for any $%
s\in \mathbb{N},$
\end{center}

\strut

where $x_{j}$'s and $y_{j}$'s are free $B$-valued random variables.

\strut

In this paper, we will consider the relation between a noncommutative
probability space and a noncommutative probability space with amalgamation
over an arbitrary unital algebra. Let $B$ be a unital algebra and let $A$ be
an algebra over $B.$ If $E:A\rightarrow B$ is the conditional expectation,
then we have the noncommutative probability space over $B,$ $(A,E).$ \ Now,
let $\varphi :A\rightarrow \mathbb{C}$ be a linear functional. Then we have
the (scalar-valued) noncommutative probability space $(A,\varphi ).$ We say
that the spaces $(A,E)$ and $(A,\varphi )$ are compatible if $\varphi
(x)=\varphi \left( E(x)\right) ,$ for all $x\in A.$\ Under compatibility, we
will consider the scalar-valued R-transform theory and amalgamated
R-transform theory.

\strut

\strut \strut

\strut

\section{Preliminaries}

\strut

\strut

\subsection{Amalgamated Free Probability Theory}

\strut

\strut

\strut \strut

In this section, we will summarize and introduced the basic results from
[19] and [20]. Throughout this section, let $B$ be a unital algebra. The
algebraic pair $(A,\varphi )$ is said to be a noncommutative probability
space over $B$ (shortly, NCPSpace over $B$) if $A$ is an algebra over $B$
(i.e $1_{B}=1_{A}\in B\subset A$) and $\varphi :A\rightarrow B$ is a $B$%
-functional (or a conditional expectation) ; $\varphi $ satisfies

\strut

\begin{center}
$\varphi (b)=b,$ for all $b\in B$
\end{center}

and

\begin{center}
$\varphi (bxb^{\prime })=b\varphi (x)b^{\prime },$ for all $b,b^{\prime }\in
B$ and $x\in A.$
\end{center}

\strut

\strut Let $(A,\varphi )$ be a NCPSpace over $B.$ Then, for the given $B$%
-functional, we can determine a moment multiplicative function $\widehat{%
\varphi }=(\varphi ^{(n)})_{n=1}^{\infty }\in I(A,B),$ where

\strut

\begin{center}
$\varphi ^{(n)}(a_{1}\otimes ...\otimes a_{n})=\varphi (a_{1}....a_{n}),$
\end{center}

\strut \strut

for all $a_{1}\otimes ...\otimes a_{n}\in A^{\otimes _{B}n},$ $\forall n\in 
\mathbb{N}.$

\strut

\strut We will denote noncrossing partitions over $\{1,...,n\}$ ($n\in 
\mathbb{N}$) by $NC(n).$ Define an ordering on $NC(n)$ ;

\strut

$\theta =\{V_{1},...,V_{k}\}\leq \pi =\{W_{1},...,W_{l}\}$\strut $\overset{%
def}{\Leftrightarrow }$ For each block $V_{j}\in \theta $, there exists only
one block $W_{p}\in \pi $ such that $V_{j}\subset W_{p},$ for $j=1,...,k$
and $p=1,...,l.$

\strut

Then $(NC(n),\leq )$ is a complete lattice with its minimal element $%
0_{n}=\{(1),...,(n)\}$ and its maximal element $1_{n}=\{(1,...,n)\}$. We
define the incidence algebra $I_{2}$ by a set of all complex-valued\
functions $\eta $ on $\cup _{n=1}^{\infty }\left( NC(n)\times NC(n)\right) $
satisfying $\eta (\theta ,\pi )=0,$ whenever $\theta \nleq \pi .$ Then,
under the convolution

\strut

\begin{center}
$*:I_{2}\times I_{2}\rightarrow \mathbb{C}$
\end{center}

defined by

\begin{center}
$\eta _{1}*\eta _{2}(\theta ,\pi )=\underset{\theta \leq \sigma \leq \pi }{%
\sum }\eta _{1}(\theta ,\sigma )\cdot \eta _{2}(\sigma ,\pi ),$
\end{center}

\strut

$I_{2}$ is indeed an algebra of complex-valued functions. Denote zeta, M\"{o}%
bius and delta functions in the incidence algebra $I_{2}$ by $\zeta ,$ $\mu $
and $\delta ,$ respectively. i.e

\strut

\begin{center}
$\zeta (\theta ,\pi )=\left\{ 
\begin{array}{lll}
1 &  & \theta \leq \pi \\ 
0 &  & otherwise,%
\end{array}
\right. $
\end{center}

\strut

\begin{center}
$\delta (\theta ,\pi )=\left\{ 
\begin{array}{lll}
1 &  & \theta =\pi \\ 
0 &  & otherwise,%
\end{array}
\right. $
\end{center}

\strut

and $\mu $ is the ($*$)-inverse of $\zeta .$ Notice that $\delta $ is the ($%
* $)-identity of $I_{2}.$ By using the same notation ($*$), we can define a
convolution between $I(A,B)$ and $I_{2}$ by

\strut

\begin{center}
$\widehat{f}\,*\,\eta \left( a_{1},...,a_{n}\,;\,\pi \right) =\underset{\pi
\in NC(n)}{\sum }\widehat{f}(\pi )(a_{1}\otimes ...\otimes a_{n})\eta (\pi
,1_{n}),$
\end{center}

\strut

where $\widehat{f}\in I(A,B)$, $\eta \in I_{2},$ $\pi \in NC(n)$ and $%
a_{j}\in A$ ($j=1,...,n$), for all $n\in \mathbb{N}.$ Notice that $\widehat{f%
}*\eta \in I(A,B),$ too. Let $\widehat{\varphi }$ be a moment multiplicative
function in $I(A,B)$ which we determined before. Then we can naturally
define a cumulant multiplicative function $\widehat{c}=(c^{(n)})_{n=1}^{%
\infty }\in I(A,B)$ by

\strut

\begin{center}
$\widehat{c}=\widehat{\varphi }*\mu $ \ \ \ or \ \ $\widehat{\varphi }=%
\widehat{c}*\zeta .$
\end{center}

\strut

\strut

This says that if we have a moment\ multiplicative function, then we always
get a cumulant multiplicative function and vice versa, by $(*).$ This
relation is so-called ''M\"{o}bius Inversion''. More precisely, we have

\strut

\begin{center}
$%
\begin{array}{ll}
\varphi (a_{1}...a_{n}) & =\varphi ^{(n)}(a_{1}\otimes ...\otimes a_{n}) \\ 
& =\underset{\pi \in NC(n)}{\sum }\widehat{c}(\pi )(a_{1}\otimes ...\otimes
a_{n})\zeta (\pi ,1_{n}) \\ 
& =\underset{\pi \in NC(n)}{\sum }\widehat{c}(\pi )(a_{1}\otimes ...\otimes
a_{n}),%
\end{array}
$
\end{center}

\strut

for all $a_{j}\in A$ and $n\in \mathbb{N}.$ Or equivalently,

\strut

\begin{center}
$%
\begin{array}{ll}
c^{(n)}(a_{1}\otimes ...\otimes a_{n}) & =\underset{\pi \in NC(n)}{\sum }%
\widehat{\varphi }(\pi )(a_{1}\otimes ...\otimes a_{n})\mu (\pi ,1_{n}).%
\end{array}
$
\end{center}

\strut \strut

Now, let $(A_{i},\varphi _{i})$ be NCPSpaces over $B,$ for all $i\in I.$
Then we can define a amalgamated free product of $A_{i}$ 's and amalgamated
free product of $\varphi _{i}$'s by

\begin{center}
$A\equiv *_{B}A_{i}$ \ \ and \ $\varphi \equiv *_{i}\varphi _{i},$
\end{center}

\strut

respectively. Then, by Voiculescu, $(A,\varphi )$ is again a NCPSpace over $%
B $ and, as a vector space, $A$ can be represented by

\begin{center}
\strut $A=B\oplus \left( \oplus _{n=1}^{\infty }\left( \underset{i_{1}\neq
...\neq i_{n}}{\oplus }(A_{i_{1}}\ominus B)\otimes ...\otimes
(A_{i_{n}}\ominus B)\right) \right) ,$
\end{center}

\strut \strut

where $A_{i_{j}}\ominus B=\ker \varphi _{i_{j}}.$ We will use Speicher's
combinatorial definition of amalgamated free product of $B$-functionals ;

\strut

\begin{definition}
Let $(A_{i},\varphi _{i})$ be NCPSpaces over $B,$ for all $i\in I.$ Then $%
\varphi =*_{i}\varphi _{i}$ is the amalgamated free product of $B$%
-functionals $\varphi _{i}$'s on $A=*_{B}A_{i}$ if the cumulant
multiplicative function $\widehat{c}=\widehat{\varphi }*\mu \in I(A,B)$ has
its restriction to $\underset{i\in I}{\cup }A_{i},$ $\underset{i\in I}{%
\oplus }\widehat{c_{i}},$ where $\widehat{c_{i}}$ is the cumulant
multiplicative function induced by $\varphi _{i},$ for all $i\in I.$ i.e

\strut

$c^{(n)}(a_{1}\otimes ...\otimes a_{n})=\left\{ 
\begin{array}{lll}
c_{i}^{(n)}(a_{1}\otimes ...\otimes a_{n}) &  & \text{if }\forall a_{j}\in
A_{i} \\ 
0_{B} &  & otherwise.%
\end{array}
\right. $
\end{definition}

\strut

Now, we will observe the freeness over $B$ ;

\strut

\begin{definition}
Let $(A,\varphi )$\strut be a NCPSpace over $B.$

\strut

(1) Subalgebras containing $B,$ $A_{i}\subset A$ ($i\in I$) are free (over $%
B $) if we let $\varphi _{i}=\varphi \mid _{A_{i}},$ for all $i\in I,$ then $%
*_{i}\varphi _{i}$ has its cumulant multiplicative function $\widehat{c}$
such that its restriction to $\underset{i\in I}{\cup }A_{i}$ is $\underset{%
i\in I}{\oplus }\widehat{c_{i}},$ where $\widehat{c_{i}}$ is the cumulant
multiplicative function induced by each $\varphi _{i},$ for all $i\in I.$

\strut

(2) Sebsets $X_{i}$ ($i\in I$) are free (over $B$) if subalgebras $A_{i}$'s
generated by $B$ and $X_{i}$'s are free in the sense of (1). i.e If we let $%
A_{i}=A\lg \left( X_{i},B\right) ,$ for all $i\in I,$ then $A_{i}$'s are
free over $B.$
\end{definition}

\strut

In [19], Speicher showed that the above combinatorial freeness with
amalgamation can be used alternatively with respect to Voiculescu's original
freeness with amalgamation.

\strut

Let $(A,\varphi )$ be a NCPSpace over $B$ and let $x_{1},...,x_{s}$ be $B$%
-valued random variables ($s\in \mathbb{N}$). Define $(i_{1},...,i_{n})$-th
moment of $x_{1},...,x_{s}$ by

\strut

\begin{center}
$\varphi (x_{i_{1}}b_{i_{2}}x_{i_{2}}...b_{i_{n}}x_{i_{n}}),$
\end{center}

\strut

for arbitrary $b_{i_{2}},...,b_{i_{n}}\in B,$ where $(i_{1},...,i_{n})\in
\{1,...,s\}^{n},$ $\forall n\in \mathbb{N}.$ Similarly, define a symmetric $%
(i_{1},...,i_{n})$-th moment by the fixed $b_{0}\in B$ by

\strut

\begin{center}
$\varphi (x_{i_{1}}b_{0}x_{i_{2}}...b_{0}x_{i_{n}}).$
\end{center}

\strut

If $b_{0}=1_{B},$ then we call this symmetric moments, trivial moments.

\strut

Cumulants defined below are main tool of combinatorial free probability
theory ; in [20], we defined the $(i_{1},...,i_{n})$-th cumulant of $%
x_{1},...,x_{s}$ by

\strut

\begin{center}
$k_{n}(x_{i_{1}},...,x_{i_{n}})=c^{(n)}(x_{i_{1}}\otimes
b_{i_{2}}x_{i_{2}}\otimes ...\otimes b_{i_{n}}x_{i_{n}}),$
\end{center}

\strut

for $b_{i_{2}},...,b_{i_{n}}\in B,$ arbitrary, and $(i_{1},...,i_{n})\in
\{1,...,s\}^{n},$ $\forall n\in \mathbb{N},$ where $\widehat{c}%
=(c^{(n)})_{n=1}^{\infty }$ is the cumulant multiplicative function induced
by $\varphi .$ Notice that, by M\"{o}bius inversion, we can always take such 
$B$-value whenever we have $(i_{1},...,i_{n})$-th moment of $%
x_{1},...,x_{s}. $ And, vice versa, if we have cumulants, then we can always
take moments. Hence we can define a symmetric $(i_{1},...,i_{n})$-th\
cumulant by $b_{0}\in B$ of $x_{1},...,x_{s}$ by

\strut

\begin{center}
$k_{n}^{symm(b_{0})}(x_{i_{1}},...,x_{i_{n}})=c^{(n)}(x_{i_{1}}\otimes
b_{0}x_{i_{2}}\otimes ...\otimes b_{0}x_{i_{n}}).$
\end{center}

\strut

If $b_{0}=1_{B},$ then it is said to be trivial cumulants of $%
x_{1},...,x_{s} $.

\strut

By Speicher, it is shown that subalgebras $A_{i}$ ($i\in I$) are free over $%
B $ if and only if all mixed cumulants vanish.

\strut

\begin{proposition}
(See [19] and [20]) Let $(A,\varphi )$ be a NCPSpace over $B$ and let $%
x_{1},...,x_{s}\in (A,\varphi )$ be $B$-valued random variables ($s\in 
\mathbb{N}$). Then $x_{1},...,x_{s}$ are free if and only if all their mixed
cumulants vanish. $\square $
\end{proposition}

\strut

\strut

\strut

\strut

\strut

\subsection{Amalgamated R-transform Theory}

\strut

\strut

In this section, we will define an R-transform of several $B$-valued random
variables. Note that to study R-transforms is to study operator-valued
distributions. R-transforms with single variable is defined by Voiculescu
(over $B,$ in particular, $B=\mathbb{C}$. See [20] and [11]). Over $\mathbb{C%
},$ Nica defined multi-variable R-transforms in [1]. In [20], we extended
his concepts, over $B.$ R-transforms of $B$-valued random variables can be
defined as $B$-formal series with its $(i_{1},...,i_{n})$-th coefficients, $%
(i_{1},...,i_{n})$-th cumulants of $B$-valued random variables, where $%
(i_{1},...,i_{n})\in \{1,...,s\}^{n},$ $\forall n\in \mathbb{N}.$

\strut

\begin{definition}
Let $(A,\varphi )$ be a NCPSpace over $B$ and let $x_{1},...,x_{s}\in
(A,\varphi )$ be $B$-valued random variables ($s\in \mathbb{N}$). Let $%
z_{1},...,z_{s}$ be noncommutative indeterminants. Define a moment series of 
$x_{1},...,x_{s}$, as a $B$-formal series, by

\strut

$M_{x_{1},...,x_{s}}(z_{1},...,z_{s})=\sum_{n=1}^{\infty }\underset{%
i_{1},..,i_{n}\in \{1,...,s\}}{\sum }\varphi
(x_{i_{1}}b_{i_{2}}x_{i_{2}}...b_{i_{n}}x_{i_{n}})\,z_{i_{1}}...z_{i_{n}},$

\strut

where $b_{i_{2}},...,b_{i_{n}}\in B$ are arbitrary for all $%
(i_{2},...,i_{n})\in \{1,...,s\}^{n-1},$ $\forall n\in \mathbb{N}.$

\strut

Define an R-transform of $x_{1},...,x_{s}$, as a $B$-formal series, by

\strut

$R_{x_{1},...,x_{s}}(z_{1},...,z_{s})=\sum_{n=1}^{\infty }\underset{%
i_{1},...,i_{n}\in \{1,...,s\}}{\sum }k_{n}(x_{i_{1}},...,x_{i_{n}})%
\,z_{i_{1}}...z_{i_{n}},$

\strut with

$k_{n}(x_{i_{1}},...,x_{i_{n}})=c^{(n)}(x_{i_{1}}\otimes
b_{i_{2}}x_{i_{2}}\otimes ...\otimes b_{i_{n}}x_{i_{n}}),$

\strut

where $b_{i_{2}},...,b_{i_{n}}\in B$ are arbitrary for all $%
(i_{2},...,i_{n})\in \{1,...,s\}^{n-1},$ $\forall n\in \mathbb{N}.$ Here, $%
\widehat{c}=(c^{(n)})_{n=1}^{\infty }$ is a cumulant multiplicative function
induced by $\varphi $ in $I(A,B).$
\end{definition}

\strut

Denote a set of all $B$-formal series with $s$-noncommutative indeterminants
($s\in \mathbb{N}$), by $\Theta _{B}^{s}$. i.e if $g\in \Theta _{B}^{s},$
then

\begin{center}
$g(z_{1},...,z_{s})=\sum_{n=1}^{\infty }\underset{i_{1},...,i_{n}\in
\{1,...,s\}}{\sum }b_{i_{1},...,i_{n}}\,z_{i_{1}}...z_{i_{n}},$
\end{center}

\strut \strut

where $b_{i_{1},...,i_{n}}\in B,$ for all $(i_{1},...,i_{n})\in
\{1,...,s\}^{n},$ $\forall n\in \mathbb{N}.$ Trivially, by definition, $%
M_{x_{1},...,x_{s}},$ $R_{x_{1},...,x_{s}}\in \Theta _{B}^{s}.$ By $\mathcal{%
R}_{B}^{s},$\ we denote a set of all R-transforms of $s$-$B$-valued random
variables. Recall that, set-theoratically,

\begin{center}
$\Theta _{B}^{s}=\mathcal{R}_{B}^{s},$ sor all $s\in \mathbb{N}.$
\end{center}

\strut \strut

By the previous section, we can also define symmetric moment series and
symmetric R-transform by $b_{0}\in B,$ by

\strut

\begin{center}
$M_{x_{1},...,x_{s}}^{symm(b_{0})}(z_{1},...,z_{s})=\sum_{n=1}^{\infty }%
\underset{i_{1},...,i_{n}\in \{1,...,s\}}{\sum }\varphi
(x_{i_{1}}b_{0}x_{i_{2}}...b_{0}x_{i_{n}})\,z_{i_{1}}...z_{i_{n}}$
\end{center}

and

\begin{center}
$R_{x_{1},...,x_{s}}^{symm(b_{0})}(z_{1},...,z_{s})=\sum_{n=1}^{\infty }%
\underset{i_{1},..,i_{n}\in \{1,...,s\}}{\sum }%
k_{n}^{symm(b_{0})}(x_{i_{1}},...,x_{i_{n}})\,z_{i_{1}}...z_{i_{n}},$
\end{center}

with

\begin{center}
$k_{n}^{symm(b_{0})}(x_{i_{1}},...,x_{i_{n}})=c^{(n)}(x_{i_{1}}\otimes
b_{0}x_{i_{2}}\otimes ...\otimes b_{0}x_{i_{n}}),$
\end{center}

for all $(i_{1},...,i_{n})\in \{1,...,s\}^{n},$ $\forall n\in \mathbb{N}.$

\strut

If $b_{0}=1_{B},$ then we have trivial moment series and trivial R-transform
of $x_{1},...,x_{s}$ denoted by $M_{x_{1},...,x_{s}}^{t}$ and $%
R_{x_{1},...,x_{s}}^{t},$ respectively.

\strut

The followings are known in [19] and [20] ;

\strut

\begin{proposition}
Let $(A,\varphi )$ be a NCPSpace over $B$ and let $%
x_{1},...,x_{s},y_{1},...,y_{p}\in (A,\varphi )$ be $B$-valued random
variables, where $s,p\in \mathbb{N}.$ Suppose that $\{x_{1},...,x_{s}\}$ and 
$\{y_{1},...,y_{p}\}$ are free in $(A,\varphi ).$ Then

\strut

(1) $%
R_{x_{1},...,x_{s},y_{1},...,y_{p}}(z_{1},...,z_{s+p})=R_{x_{1},...,x_{s}}(z_{1},...,z_{s})+R_{y_{1},...,y_{p}}(z_{s+1},...,z_{s+p}). 
$

\strut

(2) If $s=p,$ then $R_{x_{1}+y_{1},...,x_{s}+y_{s}}(z_{1},...,z_{s})=\left(
R_{x_{1},...,x_{s}}+R_{y_{1},...,y_{s}}\right) (z_{1},...,z_{s}).$

$\square $
\end{proposition}

\strut

The above proposition is proved by the characterization of freeness with
respect to cumulants. i.e $\{x_{1},...,x_{s}\}$ and $\{y_{1},...,y_{p}\}$
are free in $(A,\varphi )$ if and only if their mixed cumulants vanish. Thus
we have

\strut

$k_{n}(p_{i_{1}},...,p_{i_{n}})=c^{(n)}(p_{i_{1}}\otimes
b_{i_{2}}p_{i_{2}}\otimes ...\otimes b_{i_{n}}p_{i_{n}})$

$\ \ \ \ =\left( \widehat{c_{x}}\oplus \widehat{c_{y}}\right)
^{(n)}(p_{i_{1}}\otimes b_{i_{2}}p_{i_{2}}\otimes ...\otimes
b_{i_{n}}p_{i_{n}})$

$\ \ \ \ =\left\{ 
\begin{array}{lll}
k_{n}(x_{i_{1}},...,x_{i_{n}}) &  & or \\ 
k_{n}(y_{i_{1}},...,y_{i_{n}}) &  & 
\end{array}
\right. $

\strut

and if $s=p,$ then

$k_{n}(x_{i_{1}}+y_{i_{1}},...,x_{i_{n}}+y_{i_{n}})$

$\ =c^{(n)}\left( (x_{i_{1}}+y_{i_{1}})\otimes
b_{i_{2}}(x_{i_{2}}+y_{i_{2}})\otimes ...\otimes
b_{i_{n}}(x_{i_{n}}+y_{i_{n}})\right) $

$\ =c^{(n)}(x_{i_{1}}\otimes b_{i_{2}}x_{i_{2}}\otimes ...\otimes
b_{i_{n}}x_{i_{n}})+c^{(n)}(y_{i_{1}}\otimes b_{i_{2}}y_{i_{2}}\otimes
...\otimes b_{i_{n}}y_{i_{n}})+[Mixed]$

\strut

where $[Mixed]$ is the sum of mixed cumulants, by the bimodule map property
of $c^{(n)}$

\strut

$\ =k_{n}(x_{i_{1}},...,x_{i_{n}})+k_{n}(y_{i_{1}},...,y_{i_{n}})+0_{B}.$

\strut

Now, we will define $B$-valued boxed convolution \frame{*}, as a binary
operation on $\Theta _{B}^{s}$ ; note that if $f,g\in \Theta _{B}^{s},$ then
we can always choose free $\{x_{1},...,x_{s}\}$ and $\{y_{1},...,y_{s}\}$ in
(some) NCPSpace over $B,$ $(A,\varphi ),$ such that

\begin{center}
$f=R_{x_{1},...,x_{s}}$ \ \ and \ \ $g=R_{y_{1},...,y_{s}}.$
\end{center}

\strut

\begin{definition}
(1) Let $s\in \mathbb{N}.$ Let $(f,g)\in \Theta _{B}^{s}\times \Theta
_{B}^{s}.$ Define \frame{*}\thinspace \thinspace $:\Theta _{B}^{s}\times
\Theta _{B}^{s}\rightarrow \Theta _{B}^{s}$ by

\strut

$\left( f,g\right) =\left( R_{x_{1},...,x_{s}},\,R_{y_{1},...,y_{s}}\right)
\longmapsto R_{x_{1},...,x_{s}}\,\,\frame{*}\,\,R_{y_{1},...,y_{s}}.$

\strut

Here, $\{x_{1},...,x_{s}\}$ and $\{y_{1},...,y_{s}\}$ are free in $%
(A,\varphi )$. Suppose that

\strut

$coef_{i_{1},..,i_{n}}\left( R_{x_{1},...,x_{s}}\right)
=c^{(n)}(x_{i_{1}}\otimes b_{i_{2}}x_{i_{2}}\otimes ...\otimes
b_{i_{n}}x_{i_{n}})$

and

$coef_{i_{1},...,i_{n}}(R_{y_{1},...,y_{s}})=c^{(n)}(y_{i_{1}}\otimes
b_{i_{2}}^{\prime }y_{i_{2}}\otimes ...\otimes b_{i_{n}}^{\prime
}y_{i_{n}}), $

\strut

for all $(i_{1},...,i_{n})\in \{1,...,s\}^{n},$ $n\in \mathbb{N},$ where $%
b_{i_{j}},b_{i_{n}}^{\prime }\in B$ arbitrary. Then

\strut

$coef_{i_{1},...,i_{n}}\left( R_{x_{1},...,x_{s}}\,\,\frame{*}%
\,\,R_{y_{1},...,y_{s}}\right) $

\strut

$=\underset{\pi \in NC(n)}{\sum }\left( \widehat{c_{x}}\oplus \widehat{c_{y}}%
\right) (\pi \cup Kr(\pi ))(x_{i_{1}}\otimes y_{i_{1}}\otimes
b_{i_{2}}x_{i_{2}}\otimes b_{i_{2}}^{\prime }y_{i_{2}}\otimes ...\otimes
b_{i_{n}}x_{i_{n}}\otimes b_{i_{n}}^{\prime }y_{i_{n}})$

$\strut $

$\overset{denote}{=}\underset{\pi \in NC(n)}{\sum }\left( k_{\pi }\oplus
k_{Kr(\pi )}\right) (x_{i_{1}},y_{i_{1}},...,x_{i_{n}}y_{i_{n}}),$

\strut \strut

where $\widehat{c_{x}}\oplus \widehat{c_{y}}=\widehat{c}\mid
_{A_{x}*_{B}A_{y}},$ $A_{x}=A\lg \left( \{x_{i}\}_{i=1}^{s},B\right) $ and $%
A_{y}=A\lg \left( \{y_{i}\}_{i=1}^{s},B\right) $ and where $\pi \cup Kr(\pi
) $ is an alternating union of partitions in $NC(2n)$
\end{definition}

\strut

\begin{proposition}
(See [20])\strut Let $(A,\varphi )$ be a NCPSpace over $B$ and let $%
x_{1},...,x_{s},y_{1},...,y_{s}\in (A,\varphi )$ be $B$-valued random
variables ($s\in \mathbb{N}$). If $\{x_{1},...,x_{s}\}$ and $%
\{y_{1},...,y_{s}\} $ are free in $(A,\varphi ),$ then we have

\strut

$k_{n}(x_{i_{1}}y_{i_{1}},...,x_{i_{n}}y_{i_{n}})$

\strut

$=\underset{\pi \in NC(n)}{\sum }\left( \widehat{c_{x}}\oplus \widehat{c_{y}}%
\right) (\pi \cup Kr(\pi ))(x_{i_{1}}\otimes y_{i_{1}}\otimes
b_{i_{2}}x_{i_{2}}\otimes y_{i_{2}}\otimes ...\otimes
b_{i_{n}}x_{i_{n}}\otimes y_{i_{n}})$

\strut

$\overset{denote}{=}\underset{\pi \in NC(n)}{\sum }\left( k_{\pi }\oplus
k_{Kr(\pi )}^{symm(1_{B})}\right)
(x_{i_{1}},y_{i_{1}},...,x_{i_{n}},y_{i_{n}}),$

\strut

for all $(i_{1},...,i_{n})\in \{1,...,s\}^{n},$ $\forall n\in \mathbb{N},$ $%
b_{i_{2}},...,b_{i_{n}}\in B,$ arbitrary, where $\widehat{c_{x}}\oplus 
\widehat{c_{y}}=\widehat{c}\mid _{A_{x}*_{B}A_{y}},$ $A_{x}=A\lg \left(
\{x_{i}\}_{i=1}^{s},B\right) $ and $A_{y}=A\lg \left(
\{y_{i}\}_{i=1}^{s},B\right) .$ \ $\square $
\end{proposition}

\strut

This shows that ;

\strut

\begin{corollary}
(See [20]) Under the same condition with the previous proposition,

\strut

$R_{x_{1},...,x_{s}}\,\,\frame{*}\,%
\,R_{y_{1},...,y_{s}}^{t}=R_{x_{1}y_{1},...,x_{s}y_{s}}.$

$\square $
\end{corollary}

\strut

Notice that, in general, unless $b_{i_{2}}^{\prime }=...=b_{i_{n}}^{\prime
}=1_{B}$ in $B,$

\strut

\begin{center}
$R_{x_{1},...,x_{s}}\,\,\frame{*}\,\,R_{y_{1},...,y_{s}}\neq
R_{x_{1}y_{1},...,x_{s}y_{s}}.$
\end{center}

\strut

However, as we can see above,

\strut

\begin{center}
$R_{x_{1},...,x_{s}}\,\,\frame{*}\,%
\,R_{y_{1},...,y_{s}}^{t}=R_{x_{1}y_{1},...,x_{s}y_{s}}$
\end{center}

and

\begin{center}
$R_{x_{1},...,x_{s}}^{t}\,\,\frame{*}\,%
\,R_{y_{1},...,y_{s}}^{t}=R_{x_{1}y_{1},...,x_{s}y_{s}}^{t},$
\end{center}

\strut

where $\{x_{1},...,x_{s}\}$ and $\{y_{1},...,y_{s}\}$ are free over $B.$
Over $B=\mathbb{C},$ the last equation is proved by Nica and Speicher in [1]
and [17]. Actually, their R-transforms (over $\mathbb{C}$) is our trivial
R-transforms (over $\mathbb{C}$).\strut

\strut

\strut \strut \strut

\strut

\strut

\section{Compatiblilty of a NCPSpace and an amalgamated NCPSpace over an
algebra}

\strut \strut

\strut

\strut

In this Chapter, we will use notations defined in Section 1.3. In Section
2.1, we will introduce definitions about compatiblity. In Section 2.2, we
will observe some cumulant-relations and in Section 2.3, based on Section
2.2, we will consider the R-transform calculus. In Section 2.4, we will
observe some examples. Suppose that we have a unital algebra $B$ and an
algebra over $B,$ $A.$ Let $x_{1},...,x_{s}\in A$ be operators ($s\in 
\mathbb{N} $). Throught this paper, we will use the following notations ;
for $(i_{1},...,i_{n})\in \{1,...,s\}^{n},$ $n\in \mathbb{N},$

\strut

$\ \ k_{n}(x_{i_{1}},...,x_{i_{n}})$ \ : \ $(i_{1},...,i_{n})$-th
scalar-valued cumulants of $x_{1},...,x_{s},$

\strut

in the sense of Speicher and Nica. (In our definition, they are $%
(i_{1},...,i_{n})$-th $\mathbb{C}$-valued trivial cumulants of $%
x_{1},...,x_{s}$)

\strut

$\ \ K_{n}(x_{i_{1}},...,x_{i_{n}})$ \ : \ $(i_{1},...,i_{n})$-th $B$-valued
cumulants of $x_{1},...,x_{s}$ \ and

$\ \ K_{n}^{t}(x_{i_{1}},...,x_{i_{n}})$ \ : \ $(i_{1},...,i_{n})$-th $B$%
-valued trivial cumulants of $x_{1},...,x_{s},$

\strut

in the sense of Chapter 1.

\strut

$\ \ r_{x_{1},...,x_{s}}(z_{1},...,z_{s})$ \ : \ the scalar-valued
R-transform of $x_{1},...,x_{s},$

\strut

in the sense of Speicher and Nica. (In our definition, $%
r_{x_{1},...,x_{s}}=R_{x_{1},...,x_{s}}^{t},$ over $\mathbb{C}$)

\strut

$\ \ R_{x_{1},...,x_{s}}(z_{1},...,z_{s})$ \ : \ the $B$-valued R-transform
of $x_{1},...,x_{s}.$

\strut

Of course, $r_{x_{1},...,x_{s}}$ is an element in $\Theta _{\mathbb{C}}^{s}$
(i.e a formal series over $\mathbb{C}$) and $R_{x_{1},...,x_{s}}$ is an
element in $\Theta _{B}^{s}$ (i.e $B$-formal series). Simiralrly, we will
denote a scalar-valued moment series of $x_{1},...,x_{s}$ and a $B$-valued
moment series of $x_{1},...,x_{s}$ by

\begin{center}
\strut $m_{x_{1},...,x_{s}}$ \ \ and \ \ $M_{x_{1},...,x_{s}}$
\end{center}

respectively.

\strut

\strut \strut

\strut

\subsection{compatibility}

\strut \strut

\strut

In this section, we will introduce the compatibility. Let $B$ be a unital
algebra and let $(A,E)$ be a NCPSpace over $B.$ Also, let $(A,\varphi )$ be
a NCPSpace, in the sense of [1] and [17]. i.e an algebraic pair $(A,\varphi
) $ is a pairing of a unital algebra $A$ and a linear functional $\varphi
:A\rightarrow \mathbb{C}.$ Notice that $\varphi $ is nothing but a $\mathbb{C%
}$-functional, in the sense of Section 1.1. In [17], \'{S}niady and Speicher
introduced compatibility of $(A,\varphi )$ and $(A,E).$

\strut

\begin{definition}
Let $B$ be a unital algebra contained in an algebra $A,$ such that $%
1_{B}=1_{A},$ and let $(A,\varphi )$ be a NCPSpace. Let $(A,E)$ is a
NCPSpace over $B.$ We say $(A,\varphi )$ and $(A,E)$ are compatible if

\strut

$\varphi (a)=\varphi \left( E(a)\right) ,$ for all $a\in A.$
\end{definition}

\strut

Suppose that there exists a linear functional $\psi :B\rightarrow \mathbb{C}$
such that

\begin{center}
\strut $\varphi =\psi \circ E,$
\end{center}

then, trivially, $(A,\varphi )$ and $(A,E)$ are compatible. However, in
general, compatibility does not mean the existence of such $\psi .$

\strut

\begin{example}
Let $F_{k}=\,<g_{1},...,g_{k}>$ be a free group with $k$-generators, $%
g_{1},...,g_{k}.$ Assume that $<a,b>\,=F_{2}.$ Then we can construct a
unital group algebra $\mathbb{C}[F_{2}]$ and we can define a linear
functional $tr:\mathbb{C}[F_{2}]\rightarrow \mathbb{C},$ by

\strut

$tr\left( \underset{g\in F_{2}}{\sum }\alpha _{g}g\right) =\left\{ 
\begin{array}{lll}
\alpha _{e} &  & e\in F_{2}\text{ is the identity} \\ 
&  &  \\ 
0 &  & \text{otherwise,}%
\end{array}
\right. $

\strut

for all $\underset{g\in F_{2}}{\sum }\alpha _{g}g\in \mathbb{C}[F_{2}].$
Now, we will put

\strut

$G=\,<h\equiv aba^{-1}b^{-1}>\,\cong F_{1}=\mathbb{Z}.$

\strut

Then, similarly, we can define a unital algebra, $\mathbb{C}[F_{1}]=\mathbb{C%
}[\mathbb{Z}],$ as a subalgebra of $\mathbb{C}[F_{2}].$ Define a $\mathbb{C}%
[F_{1}]$-functional, $E:\mathbb{C}[F_{2}]\rightarrow \mathbb{C}[F_{1}]$ by

\strut

$E\left( \underset{g\in F_{2}}{\sum }\alpha _{g}g\right) =\underset{h\in
F_{1}}{\sum }\alpha _{h}h.$

\strut

Then as a NCPSpace over $\mathbb{C}[F_{1}],$ $\left( \mathbb{C}%
[F_{2}],\,\,E\right) $ is compatible with $\left( \mathbb{C}%
[F_{2}],\,tr\right) .$ Indeed,

\strut

$tr\left( E(\underset{g\in F_{2}}{\sum }\alpha _{g}g)\right) =tr\left( 
\underset{h\in F_{1}}{\sum }\alpha _{h}h\right) =\alpha _{e}=tr\left( 
\underset{g\in F_{2}}{\sum }\alpha _{g}g\right) ,$

\strut

for all $\underset{g\in F_{2}}{\sum }\alpha _{g}g\in \mathbb{C}%
[F_{2}]\hookrightarrow L(F_{2}).$
\end{example}

\strut \strut

The following lemma shows the relation between scalar-valued cumulants and a 
$B$-functional and the relation between scalar-valued moments and a $B$%
-functional. \strut In fact the follwoing lemma is just an expression gotten
from the M\"{o}bius inversion. So, at the first glance, this lemma is not so
important. However, under the compatibility, if we know $B$-moments or $B$%
-cumulants, then we can get scalar-valued moments and scalar-valued
cumulants by using them via the following lemma. For example, if we have a
NCPSpace $(A_{1}*_{B}A_{2},\,\varphi )$ and an amalgamated NCPSpace over $B,$
$(A_{1}*_{B}A_{2},\,E_{1}*E_{2})$ which are compatible and if we want to
compute scalar-valued moments (resp. scalar-valued cumulants) of operators,
then we can compute $B$-moments (resp. $B$-cumulants), first and then we can
get the scalar-valued moment series (resp. the R-transform), by using the
following lemma. In this case, it is more easy to compute $B$-moments ($B$%
-cumulants) than to compute scalar-valued moments (scalar-valued cumulants),
directly, because of the relation depending on $B,$ in $A_{1}*_{B}A_{2}.$

\strut \strut

\begin{lemma}
Let $B$ be a unital algebra and let $A$ be an algebra over $B$ (i.e $%
1_{B}=1_{A}\in B\subset A$). Let $(A,\varphi )$ be a NCPSpace and $(A,E),$ a
NCPSPace over $B,$ with a $B$-functional, $E:A\rightarrow B.$ If $(A,\varphi
)$ and $(A,E)$ are compatible and if $x_{1},...,x_{s}\in A$ are operators ($%
s\in \mathbb{N}$), then the $(i_{1},...,i_{n})$-th scalar-valued cumulant of 
$x_{1},...,x_{s}$ ($(i_{1},...,i_{n})\in \{1,...,s\}^{n},$ $n\in \mathbb{N}$%
) has the following relation with a $B$-functional $E:A\rightarrow B$ ;

\strut \strut

$k_{n}(x_{i_{1}},...,x_{i_{n}})=\underset{\pi \in NC(n)}{\sum }\left( 
\underset{(v_{1},...,v_{k})\in \pi }{\prod }\varphi \left(
E(x_{v_{1}}...x_{v_{k}})\right) \right) \mu (\pi ,1_{n}).$

\strut

Also, under the same assumption, the $(i_{1},...,i_{n})$-th moment of $%
x_{1},...,x_{s}$ has the following relation with a $B$-functional $%
E:A\rightarrow B$ ;

\strut

$\varphi (x_{i_{1}}...x_{i_{n}})=\underset{\pi \in NC(n)}{\sum }\,\underset{%
(v_{1},...,v_{k})\in \pi }{\prod }\left( \underset{\theta \in NC(k)}{\sum }\,%
\underset{(w_{1},...,w_{l})\in \theta }{\prod }\varphi \left(
E(x_{w_{1}}...x_{w_{l}})\right) \mu (\theta ,1_{k})\right) .$
\end{lemma}

\strut

\begin{proof}
By the M\"{o}bius inversion, we have that

\strut

$k_{n}(x_{i_{1}},...,x_{i_{n}})=\underset{\pi \in NC(n)}{\sum }\varphi _{\pi
}(x_{i_{1}}...x_{i_{n}})\mu (\pi ,1_{n})$

\strut

where $\varphi _{\pi }(...)$ is a partition-dependent scalar-valued moment
of $x_{1},...,x_{s},$ in the sense of Speicher and Nica

\strut

$\ \ \ \ \ \ \ =\underset{\pi \in NC(n)}{\sum }\,\left( \underset{V\in \pi }{%
\prod }\varphi _{V}(x_{i_{1}},...,x_{i_{n}})\right) \mu (\pi ,1_{n})$

\strut

where

$\varphi _{V}(x_{i_{1}},...,x_{i_{n}})\overset{def}{=}\varphi
(x_{v_{1}}...x_{v_{k}}),$

if $V=(v_{1},...,v_{k})$ is a block of $\pi $

\strut

$\ \ \ \ \ \ \ =\underset{\pi \in NC(n)}{\sum }\left( \underset{%
(v_{1},...,v_{k})\in \pi }{\prod }\varphi (x_{v_{1}}...x_{v_{k}})\right) \mu
(\pi ,1_{n})$

$\ \ \ \ \ \ \ =\underset{\pi \in NC(n)}{\sum }\left( \underset{%
(v_{1},...,v_{k})\in \pi }{\prod }\varphi \left(
E(x_{v_{1}}...x_{v_{k}})\right) \right) \mu (\pi ,1_{n}).$

\strut

Also, by the M\"{o}bius inversion, we can get that

\strut

$\varphi (x_{i_{1}}...x_{i_{n}})=\underset{\pi \in NC(n)}{\sum }k_{\pi
}(x_{i_{1}},...,x_{i_{n}})$

\strut

where $k_{\pi }(...)$ is a partition-dependent scalar-valued cumulant of $%
x_{1},...,x_{s},$ in the sense of Speicher and Nica

\strut

$\ \ \ \ =\underset{\pi \in NC(n)}{\sum }\,\underset{V\in \pi }{\prod }%
k_{V}(x_{i_{1}},...,x_{i_{n}})$

\strut

where

$k_{V}(x_{i_{1}},...,x_{i_{n}})\overset{def}{=}%
k_{k}(x_{v_{1}},...,x_{v_{k}}),$

if $V=(v_{1},...,v_{k})$ is a block of $\pi $

\strut

$\ \ \ \ =\underset{\pi \in NC(n)}{\sum }\,\underset{(v_{1},...,v_{k})\in
\pi }{\prod }k_{k}(x_{v_{1}},...,x_{v_{k}})$

$\ \ \ \ =\underset{\pi \in NC(n)}{\sum }\,\underset{(v_{1},...,v_{k})\in
\pi }{\prod }\left( \underset{\theta \in NC(k)}{\sum }\,\underset{%
(w_{1},...,w_{l})\in \theta }{\prod }\varphi \left(
E(x_{w_{1}}...x_{w_{l}})\right) \mu (\theta ,1_{k})\right) ,$

\strut

by the previous discussion for scalar-valued cumulants.
\end{proof}

\strut

\begin{example}
\strut Let $B$ be a unital algebra and $A,$ an algebra over $B.$ Suppose
that a NCPSpace $(A,\varphi )$ and an amalgamated NCPSpace over $B,$ $(A,E)$
are compatible. Let $x_{1},...,x_{7}\in A.$

\strut

(1) We can compute scalar-valued $(1,3,4)$-cumulant of them as follows ;

\strut \strut

$k_{3}(x_{1},x_{3},x_{4})=\underset{\pi \in NC(3)}{\sum }\,\left( \underset{%
(v_{1},...,v_{k})\in \pi }{\prod }\varphi \left(
E(x_{v_{1}},...,x_{v_{k}})\right) \mu (\pi ,1_{n})\right) \mu (\pi ,1_{3})$

\strut

$\ \ =\varphi \left( E(x_{1}x_{3}x_{4})\right) -\varphi \left(
E(x_{1})\right) \varphi \left( E(x_{3}x_{4})\right) -\varphi \left(
E(x_{1}x_{3})\right) \varphi \left( E(x_{4})\right) $

\strut

$\ \ \ \ \ \ \ -\varphi \left( E(x_{1}x_{4})\right) \varphi \left(
E(x_{3})\right) +\varphi \left( E(x_{1})\right) \varphi \left(
E(x_{3})\right) \varphi \left( E(x_{4})\right) .$

\strut

(2) Now, suppose that each $x_{j}$, \ $j=1,...,7,$ is centered over $B.$
(i.e $E(x_{j})=0_{B},$ for all $j=1,...,7.$ ) Then

\strut

$k_{3}(x_{1},x_{3},x_{4})=\varphi \left( E(x_{1}x_{3}x_{4})\right) $

\strut

$\ \ \ \ \ \ =\varphi \left( \underset{\theta \in NC(3)}{\sum }\widehat{C}%
(\theta )(x_{1}\otimes x_{3}\otimes x_{4})\right) $

\strut

$\ \ \ \ \ \ =\varphi \left( K_{3}(x_{1},x_{3},x_{4})\right)
+K_{1}(x_{1})K_{2}(x_{3},x_{4})+K_{2}\left( x_{1},K_{1}(x_{3})x_{4}\right) $

$\ \ \ \ \ \ \ \ \ \ \ \ \ \ \ \ \ \ \ \ \ \ \ \ \ \ \ \ \ \ \ \ \ \ \ \ \
+K_{2}(x_{1},x_{3})K_{1}(x_{4})+K_{1}(x_{1})K_{1}(x_{2})K_{1}(x_{3})$

$\ \ \ \ \ \ =\varphi \left( K_{3}(x_{1},x_{3},x_{4})\right) ,$

\strut

since $E(x_{j})=K_{1}(x_{j}),$ for all $j=1,...,7.$ So, if $x_{j}$'s are
centered, then

\strut

$k_{3}(x_{1},x_{3},x_{4})=\varphi \left( K_{3}(x_{1},x_{3},x_{4})\right) .$
\end{example}

\strut

We hope that scalar-valued cumulants and operator-valued cumulants have nice
property under the compatibility such as

\strut

\begin{center}
$\varphi \left( K_{n}^{t}(x_{i_{1}},...,x_{i_{n}})\right)
=k_{n}(x_{i_{1}},...,x_{i_{n}}).$
\end{center}

\strut

However, in general, the above relation does NOT hold true (See the above
example (1)). Observe the following ; if $x,y\in A$ are operators, then

\strut

\begin{center}
$%
\begin{array}{ll}
\varphi \left( K_{2}^{t}(x,y)\right) & =\varphi \left( E(xy)-E(x)E(y)\right)
\\ 
& =\varphi \left( E(xy)\right) -\varphi \left( E(x)E(y)\right)%
\end{array}
$
\end{center}

and

\begin{center}
$%
\begin{array}{ll}
k_{2}(xy) & =\varphi (xy)-\varphi (x)\varphi (y) \\ 
& =\varphi \left( E(xy)\right) -\varphi \left( E(x)\right) \varphi \left(
E(y)\right) .%
\end{array}
$
\end{center}

\strut

Thus to get $\varphi \left( K_{2}^{t}(x,y)\right) =k_{2}(x,y),$ we need the
following eaulity ;

\strut

\begin{center}
$\varphi \left( E(x)E(y)\right) =\varphi \left( E(x)\right) \varphi \left(
E(y)\right) ,$
\end{center}

\strut

since

$\ \ \ \ \ \ \ \varphi \left( K_{2}^{t}(x,y)\right) -k_{2}(x,y)$

$\ \ \ \ \ \ \ \ \ =\varphi \left( E(xy)\right) -\varphi \left(
E(x)E(y)\right) -\left( \varphi \left( E(xy)\right) -\varphi \left(
E(x)\right) \varphi \left( E(y)\right) \right) $

$\ \ \ \ \ \ \ \ \ =\varphi \left( E(x)\right) \varphi \left( E(y)\right)
-\varphi \left( E(x)E(y)\right) .$

\strut

If $\varphi :A\rightarrow \mathbb{C}$ is a homomorphism, then it happens,
however, in general, we cannot guarantee the above equality. For example,
suppose that we have a UHF-algebra $A$ and $B=M_{N}(\mathbb{C})$ \ ($N\neq 1$%
) and assume that $E_{N}:A\rightarrow M_{N}(\mathbb{C})$ is a canonical
conditional expectation and $\varphi =\underset{k\rightarrow \infty }{\lim }%
\varphi _{k}:A\rightarrow \mathbb{C}$ is a normalized trace. Let

\begin{center}
$E_{N}(x)=\left( 
\begin{array}{lllll}
0 &  &  &  & O \\ 
1 & 0 &  &  &  \\ 
0 & \ddots & 0 &  &  \\ 
& \ddots & \ddots & \ddots &  \\ 
O &  & 0 & 1 & 0%
\end{array}
\right) $ \ and \ $E_{N}(y)=\left( 
\begin{array}{lllll}
0 & 1 & 0 &  & O \\ 
& 0 & \ddots & \ddots &  \\ 
&  & 0 & \ddots & 0 \\ 
&  &  & \ddots & 1 \\ 
O &  &  &  & 0%
\end{array}
\right) $
\end{center}

\strut

in $M_{N}(\mathbb{C}).$ Then

\strut

\begin{center}
$\varphi \left( E(x)E(y)\right) =\varphi _{N}\left( \left( 
\begin{array}{llll}
0 &  &  & O \\ 
& 1 &  &  \\ 
&  & \ddots &  \\ 
O &  &  & 1%
\end{array}
\right) \right) =\frac{N-1}{N}$
\end{center}

and

\begin{center}
$\varphi \left( E(x)\right) \varphi \left( E(y)\right) =\varphi _{N}\left(
E(x)\right) \cdot \varphi _{N}\left( E(y)\right) =0.$
\end{center}

\strut

More generally, to satisfy

\strut

\begin{center}
$\varphi \left( K_{n}^{t}(x_{i_{1}},...,x_{i_{n}})\right)
=k_{n}(x_{i_{1}},...,x_{i_{n}}),$
\end{center}

\strut

operators $x_{1},...,x_{s}\in A$ should satisfy ;

\strut

(*) : $\varphi \left( \underset{V\in \pi (o)}{\prod }\widehat{E}(\pi \mid
_{V})(x_{i_{1}}\otimes ...\otimes \widehat{E}(\pi \mid
_{W})(x_{i_{1}}\otimes ...\otimes x_{i_{n}})x_{i_{j}}\otimes ...\otimes
x_{i_{n}})\right) $

$=\underset{V\in \pi (o)}{\prod }\varphi \left( \widehat{E}(\pi \mid
_{V})(x_{i_{1}}\otimes ...\otimes \varphi \left( \widehat{E}(\pi \mid
_{W})(x_{i_{1}}\otimes ...\otimes x_{i_{n}})\right) x_{i_{j}}\otimes
...\otimes x_{i_{n}})\right) ,$

\strut

for any $\pi \in NC(n),$ where $W$ is an inner block having its outer block $%
V,$ for the fixed $\pi \in NC(n).$

\strut

\begin{proposition}
Let $B$ be a unital algebra and $A,$ an algebra over $B.$ Suppose that a
NCPSpace $(A,\varphi )$ and an amalgamated NCPSpace over $B,$ $(A,E)$ are
compatible. Assume that operators $x_{1},...,x_{s}\in A$ ($s\in \mathbb{N}$%
)\ satisfy

\strut

(*) : $\varphi \left( \underset{V\in \pi (o)}{\prod }\widehat{E}(\pi \mid
_{V})(x_{i_{1}}\otimes ...\otimes \widehat{E}(\pi \mid
_{W})(x_{i_{1}}\otimes ...\otimes x_{i_{n}})x_{i_{j}}\otimes ...\otimes
x_{i_{n}})\right) $

$=\underset{V\in \pi (o)}{\prod }\varphi \left( \widehat{E}(\pi \mid
_{V})(x_{i_{1}}\otimes ...\otimes \varphi \left( \widehat{E}(\pi \mid
_{W})(x_{i_{1}}\otimes ...\otimes x_{i_{n}})\right) x_{i_{j}}\otimes
...\otimes x_{i_{n}})\right) ,$

\strut

for any $\pi \in NC(n),$ for all $(i_{1},...,i_{n})\in \{1,...,s\}^{n},$ $%
n\in \mathbb{N}$, where $W\in \pi (i)$ is an inner block of $V\in \pi (o),$
for the fixed $\pi \in NC(n).$ Then

\strut \strut

$\varphi \left( K_{n}^{t}(x_{i_{1}},...,x_{i_{n}})\right)
=k_{n}(x_{i_{1}},...,x_{i_{n}}),$

\strut

for all $(i_{1},...,i_{n})\in \{1,...,s\}^{n},$ $n\in \mathbb{N}.$ In
particular,

\strut

$r_{x_{1},...,x_{s}}(z_{1},...,z_{s})=\sum_{n=1}^{\infty }\underset{%
i_{1},...,i_{n}\in \{1,...,s\}}{\sum }\varphi \left(
coef_{i_{1},...,i_{n}}(R_{x_{1},...,x_{s}}^{t})\right)
z_{i_{1}}...z_{i_{n}}. $

$\square $
\end{proposition}

\strut

\begin{quote}
\strut \frame{\textbf{Notation}} We will say that operators $%
x_{1},...,x_{s}\in A$ satisfy property (*) if $x_{1},...,x_{s}$ satisfy the
relation (*) introduced in the previous proposition.

\strut
\end{quote}

\begin{corollary}
Let $B$ be a unital algebra and $A,$ an algebra over $B.$ Suppose that a
NCPSpace $(A,\varphi )$ and an amalgamated NCPSpace over $B,$ $(A,E)$ are
compatible. If a linear functional $\varphi :A\rightarrow \mathbb{C}$ is an
algebra homomorphism, then, for $x_{1},...,x_{s}\in A$ ($s\in \mathbb{N}$),

\strut

$r_{x_{1},...,x_{s}}(z_{1},...,z_{s})=\sum_{n=1}^{\infty }\underset{%
i_{1},...,i_{n}\in \{1,...,s\}}{\sum }\varphi \left(
coef_{i_{1},...,i_{n}}(R_{x_{1},...,x_{s}}^{t})\right)
z_{i_{1}}...z_{i_{n}}. $

$\square $
\end{corollary}

\strut \strut

This shows that it is difficult to verify the relation between $%
r_{x_{1},...,x_{s}}$ and $R_{x_{1},...,x_{s}}.$ This also says that $B$%
-freeness and scalar-valued freeness have a deep gap, in general.

\strut

\begin{theorem}
Let $B$ be a unital algebra and $A,$ an algebra over $B.$ Suppose that a
NCPSpace $(A,\varphi )$ and an amalgamated NCPSpace over $B,$ $(A,E)$ are
compatible. If $\{x_{1},...,x_{s}\}$ and $\{y_{1},...,y_{t}\}$ are $B$-free
families of $B$-valued random variables ($s,t\in \mathbb{N}$) and if $%
x_{1},...,x_{s},y_{1},...,y_{t}$ satisfy property(*), then $%
\{x_{1},...,x_{s}\}$ and $\{y_{1},...,y_{t}\}$ are free in $(A,\varphi ).$
In particular, if a linear functional $\varphi :A\rightarrow \mathbb{C}$ is
nondegenerated and if $\{x_{1},...,x_{s}\}$ and $\{y_{1},...,y_{t}\}$ are
central in the sense of [20], then the converse also holds true.
\end{theorem}

\strut

\begin{proof}
Assume that operators $x_{1},...,x_{s},y_{1},...,y_{t}$ satisfies property
(*). Then, by the previous proposition, we have that

\strut

$k_{n}(p_{1},...,p_{n})=\varphi \left( K_{n}^{t}(p_{1},...,p_{n})\right) ,$

\strut

where $p_{1},...,p_{n}\in \{x_{1},...,x_{s}\}\cup \{y_{1},...,y_{t}\}.$
Thus, if $\{x_{1},...,x_{s}\}$ and $\{y_{1},...,y_{t}\}$ are $B$-free over $%
B,$ in $(A,E),$ then a mixed cumulant of $x_{1},...,x_{s},y_{1},...,y_{t}$
vanishs ;

\strut

$k_{n}(p_{1},...,p_{n})=\varphi (0_{B})=0.$

\strut

This shows that all mixed scalar-valued cumulants of $%
x_{1},...,x_{s},y_{1},...,y_{t}$ vanish, too. Equivalently, $%
\{x_{1},...,x_{s}\}$ and $\{y_{1},...,y_{t}\}$ are free in $(A,\varphi ).$

\strut \strut

Now, suppose that the linear functional $\varphi $ is nondegenerated and
assume that $x_{1},...,x_{s}$ and $y_{1},...,y_{t}$ are central. By the
property ''central'', $B$-freeness is equivalent to the statement [all mixed
trivial cumulants vanish] (See [20] and [29]). So, we have that ; for the
mixed case,

\strut

$0=k_{n}(p_{1},...,p_{n})=\varphi \left( K_{n}^{t}(p_{1},...,p_{n})\right) $

$\Longleftrightarrow $

$K_{n}^{t}(p_{1},...,p_{n})=0_{B},$

\strut

by the nondegeneratedness of $\varphi .$ Therefore, $\{x_{1},...,x_{s}\}$
and $\{y_{1},...,y_{t}\}$ are $B$-free over $B$ in $(A,E)$
\end{proof}

\strut

Also, the above theorem says that even under the compatibility, $B$-freeness
and scalar-valued freeness have a deep gap, in general. Actually, the
property (*) is a very powerful assumption. We redefine the property (*)
which is the extended concept of having the property (*)

\strut

\begin{definition}
Let $B$ be a unital algebra and $A,$ an algebra over $B.$ Suppose that a
NCPSpace $(A,\varphi )$ and an amalgamated NCPSpace over $B,$ $(A,E)$ are
compatible. We say that operators $x_{1},...,x_{s}\in A$ satisfy the
property (*), if

\strut

$k_{n}(x_{i_{1}},...,x_{i_{n}})=\varphi \left(
K_{n}^{t}(x_{i_{1}},...,x_{i_{n}})\right) ,$

\strut

for all $(i_{1},...,i_{n})\in \{1,...,s\}^{n},$ $n\in \mathbb{N}.$
\end{definition}

\strut

\strut

\strut

\strut

\subsection{R-transformTheory under the compitablity : Connection between
scalar-valued R-transforms and operator-valued Moment Series}

\strut

\strut

\strut

In this section, we will observe the connections between scalar-valued
R-transforms and operator-valued moment series under the compatibility of a
NCPSpace and amalgamated NCPSpace. This connection will be used for studying
scalar-valued distributions of (some) operators and operator-valued
distribution of those operators. From the previous section, we can compute
operator-valued R-transforms and scalar-valued R-transforms. Let $B$ be a
unital algebra.

\strut \strut

\begin{theorem}
(See Theorem 6, in [17]) Let $B$ be a unital algebra and $A,$ an algebra
over $B.$ Let $(A,\varphi )$ and $(A,E)$ be compatible and let $%
x_{1},...,x_{s}\in A$ be operators ($s\in \mathbb{N}$). Then $%
\{x_{1},...,x_{s}\}$ and $B$ are free in $(A,\varphi )$ if and only if we
have that

\strut

$K_{n}(x_{i_{1}},...,x_{i_{n}})=\left( \varphi (b_{i_{2}})\cdot \cdot \cdot
\varphi (b_{i_{n}})\right) k_{n}(x_{i_{1}},...,x_{i_{n}})\,\,\cdot
\,1_{B}\in B,$

\strut

where

\strut

$K_{n}(x_{i_{1}},...,x_{i_{n}})=C^{(n)}\left( x_{i_{1}}\otimes
b_{i_{2}}x_{i_{2}}\otimes ...\otimes b_{i_{n}}x_{i_{n}}\right) ,$

\strut

for arbitrary $b_{i_{2}},...,b_{i_{n}}\in B,$ and for all $%
(i_{1},...,i_{n})\in \{1,...,s\}^{n},$ $n\in \mathbb{N}.$ \ $\square $
\end{theorem}

\strut \strut

The above theorem shows that if some family of operators in $A$ and a unital
algebra $B,$ in $A,$ are free in $(A,\varphi ),$ then an operator-valued ($B$%
-valued) cumulants of those operators are easily expressed in termes of
multiplication of some scalar-value and scalar-valued trivial cumulants of
those operators.

\strut

\begin{corollary}
Let $B$ be a unital algebra and $A,$ an algebra over $B.$ Let $(A,\varphi )$
and $(A,E)$ be compatible NCPSpace and NCPSpace over $B.$ If $%
x_{1},...,x_{s}\in A$ are operators ($s\in \mathbb{N}$) and if $%
\{x_{1},...,x_{s}\}$ and $B$ are free in $(A,\varphi ),$ then there exist $%
\alpha _{i_{1},...,i_{n}}\in \mathbb{C},$ for all $(i_{1},...,i_{n})\in
\{1,...,s\}^{n},$ $n\in \mathbb{N},$ such that

\strut

$R_{x_{1},...,x_{s}}(z_{1},...,z_{s})=\sum_{n=1}^{\infty }\underset{%
i_{1},...,i_{n}\in \{1,...,s\}}{\sum }\left( \alpha
_{i_{1},...,i_{n}}\,k_{n}(x_{i_{1}},...,x_{i_{n}})\cdot 1_{B}\right)
z_{i_{1}}...z_{i_{n}},$

\strut

in $\Theta _{B}^{s},$ as a $B$-formal series. \ $\square $
\end{corollary}

\strut

\strut By the previous theorem, in particular, we can characterize $\alpha
_{i_{1},...,i_{n}}\in \mathbb{C},$ for each $(i_{1},...,i_{n})\in
\{1,...,s\}^{n}\dot{,}$ $n\in \mathbb{N},$ as follows ;

\strut

\begin{center}
$\alpha _{i_{1},...,i_{n}}=\varphi (1_{B})\cdot \varphi (b_{i_{2}})\cdot
\cdot \cdot \varphi (b_{i_{n}})\in \mathbb{C},$
\end{center}

\strut

where $b_{i_{2}},...,b_{i_{n}}\in B\subset A$ are determined in

\strut

\begin{center}
$K_{n}(x_{i_{1}},...,x_{i_{n}})=C^{(n)}\left( x_{i_{1}}\otimes
b_{i_{2}}x_{i_{2}}\otimes ...\otimes b_{i_{n}}x_{i_{n}}\right) .$
\end{center}

\strut \strut

Notice that the above $B$-valued R-transform of $x_{1},...,x_{s},$ $%
R_{x_{1},...,x_{s}}$ is a $B$-formal series, but each coefficient has the
form of $\alpha \cdot 1_{B},$ where $\alpha \in \mathbb{C}.$ This says that
we can regard the R-transform, $R_{x_{1},...,x_{s}}$, as a scalar-valued
formal series, with its $(i_{1},...,i_{n})$-th coefficients

\strut

\begin{center}
$k_{n}(x_{i_{1}},\,\varphi (b_{i_{2}})x_{i_{2}},...,\varphi
(b_{i_{n}})x_{i_{n}})=\left( \varphi (b_{i_{2}})...\varphi
(b_{i_{n}})\right) k_{n}(x_{i_{1}},...,x_{i_{n}}).$
\end{center}

\strut \strut

Now, we will denote boxed convolution and amalgamated boxed convolution over
an algebra $B$ by \ \strut \frame{*} \ \ \ and \ \ \ \frame{*}$_{B},$
respectively.

\strut

\strut

\begin{theorem}
Let $B$ be a unital algebra and $A,$ an algebra over $B.$ Let $(A,\varphi )$
and $(A,E)$ be compatible. If $x_{1},...,x_{s}\in A$ are operators ($s\in 
\mathbb{N}$) and if $\{x_{1},...,x_{s}\}$ are free from $B,$ in $(A,\varphi
),$ then there exists a $B$-formal series $g\in \Theta _{B}^{s}$ such that

\strut

$r_{x_{1},...,x_{s}}(z_{1},...,z_{s})=\sum_{n=1}^{\infty }\underset{%
i_{1},...,i_{n}\in \{1,...,s\}}{\sum }[$scalar part of $%
K_{n}^{t}(x_{i_{1}},...,x_{i_{n}})]z_{i_{1}}...z_{i_{n}}.$
\end{theorem}

\strut

\begin{proof}
By the previous theorem, we have that, for $(i_{1},...,i_{n})\in
\{1,...,s\}^{n},$ $n\in \mathbb{N},$

\strut

$K_{n}^{t}(x_{i_{1}},...,x_{i_{n}})=k_{n}(x_{i_{1}},...,x_{i_{n}})\cdot
1_{B}\equiv c_{i_{1},...,i_{n}}\cdot 1_{B}.$

\strut

Therefore, in this case, we can get that

\strut

$r_{x_{1},...,x_{s}}(z_{1},...,z_{s})=\sum_{n=1}^{\infty }\underset{%
i_{1},...,i_{n}\in \{1,...,s\}}{\sum }%
k_{n}(x_{i_{1}},...,x_{i_{n}})z_{i_{1}}...z_{i_{n}}$

$\ \ =\sum_{n=1}^{\infty }\underset{i_{1},...,i_{n}\in \{1,...,s\}}{\sum }%
c_{i_{1},...,i_{n}}\,z_{i_{1}}...z_{i_{n}}$

$\ \ =\sum_{n=1}^{\infty }\underset{i_{1},...,i_{n}\in \{1,...,s\}}{\sum }[$%
scalar part of $K_{n}^{t}(x_{i_{1}},...,x_{i_{n}})]z_{i_{1}}...z_{i_{n}}$
\end{proof}

\strut \strut

\begin{remark}
Remark that if a set of operators$\ \{x_{1},...,x_{s}\}\subset A$ is free
from $B,$ in $(A,\varphi ),$ then $x_{1},...,x_{s}$ satisfy the property
(*). Indeed,

\strut

$%
\begin{array}{ll}
\varphi \left( K_{n}^{t}(x_{i_{1}},...,x_{i_{n}})\right) & =\varphi \left(
k_{n}(x_{i_{1}},...,x_{i_{n}})\cdot 1_{B}\right) \\ 
& =k_{n}(x_{i_{1}},...,x_{i_{n}})\cdot \varphi (1_{B}) \\ 
& =k_{n}(x_{i_{1}},...,x_{i_{n}}),%
\end{array}
$

\strut

for all $(i_{1},...,i_{n})\in \{1,...,s\}^{n},$ $n\in \mathbb{N}.$
\end{remark}

\strut

\begin{corollary}
Let $B$ be a unital algebra and $A,$ an algebra over $B.$ Suppose a NCPSpace 
$(A,\varphi )$ and an amalgamated NCPSpace over $B,$ $(A,E)$ are compatible.
Let $x_{1},...,x_{s},y_{1},...,y_{s}\in A$ be operators such that $%
\{x_{1},...,x_{s}\}$ and $\{y_{1},...,y_{s}\}$ are $B$-free in $(A,E).$ If $%
\{x_{1},...,x_{s},y_{1},...,y_{s}\}$ and $B$ are free in $(A,\varphi ),$ then

\strut

(1) $%
r_{x_{1},...,x_{s},y_{1},...,y_{s}}(z_{1},...,z_{2s})=r_{x_{1},..,x_{s}}(z_{1},...,z_{s})+r_{y_{1},...,y_{s}}(z_{s+1},...,z_{2s}). 
$

\strut

(2) $%
r_{x_{1}+y_{1},...,x_{s}+y_{s}}(z_{1},...,z_{s})=r_{x_{1},...,x_{s}}(z_{1},...,z_{s})+r_{y_{1},...,y_{s}}(z_{1},...,z_{s}) 
$

\strut

(3) $r_{x_{1}y_{1},...,x_{s}y_{s}}(z_{1},...,z_{s})=\left(
r_{x_{1},...,x_{s}}\,\,\frame{*}\,\,r_{y_{1},...,y_{s}}\right)
(z_{1},...,z_{s}).$
\end{corollary}

\strut

\begin{proof}
(1) Since a set of operators in $A,$ $\{x_{1},...,x_{s},y_{1},...,y_{s}\}$
are free from $B,$ in $(A,\varphi ),$ by the previous theorem, we can get
that

\strut

$r_{x_{1},...,x_{s},y_{1},...,y_{s}}(z_{1},...,z_{2s})=\sum_{n=1}^{\infty }%
\underset{i_{1},...,i_{n}\in \{1,...,2s\}}{\sum }c_{i_{1},...,i_{n}}%
\,z_{i_{1}}...z_{i_{n}}$

\strut

where

$c_{i_{1},...,i_{n}}=$ the scalar part of $%
K_{n}^{t}(p_{i_{1}},...,p_{i_{n}}) $

with

$K_{n}^{t}(p_{i_{1}},...,p_{i_{n}})=\left\{ 
\begin{array}{llll}
K_{n}^{t}(x_{i_{1}},...,x_{i_{n}}) &  &  & \text{or} \\ 
&  &  &  \\ 
K_{n}^{t}(y_{i_{1}},...,y_{i_{n}}), &  &  & 
\end{array}
\right. $

\strut

by the $B$-freeness of $\{x_{1},...,x_{s}\}$ and $\{y_{1},...,y_{s}\}.$ Put

\strut

$K_{n}^{t}(x_{i_{1}},...,x_{i_{n}})=c_{i_{1},...,i_{n}}^{x}\cdot 1_{B}$ \
and \ $K_{n}^{t}(y_{i_{1}},...,y_{i_{n}})=c_{i_{1},...,i_{n}}^{y}\cdot
1_{B}, $

\strut

for all $(i_{1},...,i_{n})\in \{1,...,2s\}^{n},$ $n\in \mathbb{N}.$

\strut

Thus, we have that

\strut

$r_{x_{1},...,x_{s},y_{1},...,y_{s}}(z_{1},...,z_{2s})$

$\ \ \ \ \ \ \ \ \ =\sum_{n=1}^{\infty }\underset{i_{1},...,i_{n}\in
\{1,...,2s\}}{\sum }\left(
c_{i_{1},...,i_{n}}^{x}+c_{i_{1},...,i_{n}}^{y}\right) z_{i_{1}}...z_{i_{n}}$

$\ \ \ \ \ \ \ \ \
=r_{x_{1},...,x_{s}}(z_{1},...,z_{s})+r_{y_{1},...,y_{s}}(z_{s+1},...,z_{2s}). 
$

\strut

(2) Since $\{x_{1},...,x_{s},y_{1},...,y_{s}\}$ are free from $B,$ in $%
(A,\varphi ),$ so are $\{x_{1}+y_{1},...,x_{s}+y_{s}\}$ and $B.$ Thus

\strut

$r_{x_{1}+y_{1},...,x_{s}+y_{s}}(z_{1},...,z_{s})=\sum_{n=1}^{\infty }%
\underset{i_{1},...,i_{n}\in \{1,...,s\}}{\sum }%
c_{i_{1},...,i_{n}}z_{i_{1}}...z_{i_{n}}$

\strut

where

$c_{i_{1},...,i_{n}}=$ the scalar part of $%
K_{n}^{t}(x_{i_{1}}+y_{i_{1}},...,x_{i_{n}}+y_{i_{n}}).$

\strut

By the $B$-freeness of $\{x_{1},...,x_{s}\}$ and $\{y_{1},...,y_{s}\},$

\strut

$%
K_{n}^{t}(x_{i_{1}}+y_{i_{1}},...,x_{i_{n}}+y_{i_{n}})=K_{n}^{t}(x_{i_{1}},...,x_{i_{n}})+K_{n}^{t}(y_{i_{1}},...,y_{i_{n}}), 
$

\strut

for all $(i_{1},...,i_{n})\in \{1,...,s\}^{n},$ $n\in \mathbb{N}.$ Notice
that both $\{x_{1},...,x_{s}\}$ and $\{y_{1},...,y_{s}\}$ are free from $B,$
in $(A,\varphi )$. Hence

\strut

$K_{n}^{t}(x_{i_{1}},...,x_{i_{n}})=k_{n}(x_{i_{1}},...,x_{i_{n}})\cdot
1_{B}\equiv c_{i_{1},...,i_{n}}^{x}\cdot 1_{B}$

and

$K_{n}^{t}(y_{i_{1}},...,y_{i_{n}})=k_{n}(y_{i_{1}},...,y_{i_{n}})\cdot
1_{B}\equiv c_{i_{1},...,i_{n}}^{y}\cdot 1_{B}.$

\strut

Therefore,

\strut

$c_{i_{1},...,i_{n}}\cdot 1_{B}=\left(
c_{i_{1},...,i_{n}}^{x}+c_{i_{1},...,i_{n}}^{y}\right) \cdot 1_{B}$

and\strut

$r_{x_{1}+y_{1},...,x_{s}+y_{s}}(z_{1},...,z_{s})=\left(
r_{x_{1},...,x_{s}}+r_{y_{1},...,y_{s}}\right) (z_{1},....,z_{s}).$

\strut

(3) Similarly, $\{x_{1}y_{1},...,x_{s}y_{s}\}$ and $B$ are free in $%
(A,\varphi ).$ Thus

\strut

$r_{x_{1}y_{1},...,x_{s}y_{s}}(z_{1},...,z_{s})=\sum_{n=1}^{\infty }\underset%
{i_{1},...,i_{n}\in \{1,...,s\}}{\sum }%
c_{i_{1},...,i_{n}}z_{i_{1}}...z_{i_{n}}$

with

$c_{i_{1},...,i_{n}}=$ the scalar part of $%
K_{n}^{t}(x_{i_{1}}y_{i_{1}},...,x_{i_{n}}y_{i_{n}}).$

\strut \strut

By the $B$-freeness of $\{x_{1},...,x_{s}\}$ and $\{y_{1},...,y_{s}\},$ we
have that

\strut

$K_{n}^{t}(x_{i_{1}}y_{i_{1}},...,x_{i_{n}}y_{i_{n}})=C^{(n)}\left(
x_{i_{1}}y_{i_{1}}\otimes ...\otimes x_{i_{n}}y_{i_{n}}\right) $

\strut

$\ \ \ \ =\underset{\pi \in NC(n)}{\sum }\left( \widehat{C_{x}}\oplus 
\widehat{C_{y}}\right) (\pi \cup Kr(\pi ))(x_{i_{1}}\otimes y_{i_{1}}\otimes
...\otimes x_{i_{n}}\otimes y_{i_{n}})$

\strut

$\ \ \ \ =coef_{i_{1},...,i_{n}}\left( R_{x_{1},...,x_{s}}^{t}\,\,\frame{*}%
_{B}\,\,R_{y_{1},...,y_{s}}^{t}\right) ,$

\strut

where $\widehat{C_{x}}\oplus \widehat{C_{y}}=\widehat{C}\mid
_{A_{x}*_{B}A_{y}},$ $A_{x}=A\lg \left( \{x_{j}\}_{j=1}^{s},B\right) $ and $%
A_{y}=A\lg \left( \{y_{j}\}_{j=1}^{s},B\right) .$ Notice that, since both $%
\{x_{1},...,x_{s}\}$ and $\{y_{1},...,y_{s}\}$ are free from $B,$ in $%
(A,\varphi ),$

\strut

$K_{N}^{t}\left( x_{j_{1}},...,x_{j_{N}}\right) =k_{N}\left(
x_{j_{1}},...,x_{j_{N}}\right) \cdot 1_{B}$

and

$K_{N}^{t}\left( y_{j_{1}},...,y_{j_{N}}\right) =k_{N}\left(
y_{j_{1}},...,y_{j_{N}}\right) \cdot 1_{B}.$

\strut

(i.e they are multiplications of scalar-value and $1_{B}$ and hence $%
x_{1},...,x_{s}$ and $y_{1},...,y_{s}$ are central !) Thus we can get that

\strut

$\underset{\pi \in NC(n)}{\sum }\left( \widehat{C_{x}}\oplus \widehat{C_{y}}%
\right) (\pi \cup Kr(\pi ))(x_{i_{1}}\otimes y_{i_{1}}\otimes ...\otimes
x_{i_{n}}\otimes y_{i_{n}})$

\strut

$\ \ \ \ \ \ =\underset{\pi \in NC(n)}{\sum }\left( \widehat{C_{x}}(\pi
)(x_{i_{1}}\otimes ...\otimes x_{i_{n}})\right) \left( \widehat{C_{y}}%
(Kr(\pi ))(y_{i_{1}}\otimes ...\otimes y_{i_{n}})\right) $

$\ \ \ \ \ \ =\underset{\pi \in NC(n)}{\sum }\left( k_{\pi
}(x_{i_{1}},...,x_{i_{n}})\cdot 1_{B}\right) \left( k_{Kr(\pi
)}(y_{i_{1}},...,y_{i_{n}})\cdot 1_{B}\right) .$

\strut

Therefore,

\strut

$r_{x_{1}y_{1},...,x_{s}y_{s}}(z_{1},...,z_{s})=\left(
r_{x_{1},...,x_{s}}\,\,\frame{*}\,\,r_{y_{1},...,y_{s}}\right)
(z_{1},...,z_{s})$
\end{proof}

\strut \strut

From now, we will consider the one-variable scalar-valued R-transforms and
operator-valued ones. As we have seen, if we have compatible $(A,\varphi )$
and $(A,E)$, over $B$ and if $x_{1},...,x_{s}$ are operators ($s\in \mathbb{N%
}$) in $A,$ then $(i_{1},...,i_{n})$-th scalar-valued cumulants of $%
x_{1},...,x_{s}$ are

\strut

\begin{center}
$k_{n}(x_{i_{1}},...,x_{i_{n}})=\underset{\pi \in NC(n)}{\sum }\,\underset{%
(v_{1},...,v_{k})\in \pi }{\prod }\varphi \left(
E(x_{v_{1}}...x_{v_{k}})\right) \mu (\pi ,1_{n}).$
\end{center}

\strut

Now, let $x\in A$ be an operator. Then we can consider the $n$-th cumulants
of $x,$ for all $n\in \mathbb{N}$ and, under the compatibility, we can get
that

\strut

\begin{center}
$k_{n}\left( \underset{n-times}{\underbrace{x,.....,x}}\right) =\underset{%
\pi \in NC(n)}{\sum }\,\underset{(v_{1},...,v_{k})\in \pi }{\prod }\varphi
\left( E(x^{k})\right) \mu (\pi ,1_{n}).$
\end{center}

\strut

Later, what we are most interested in is to compute $n$-th moments of an
operator $x$ in $A.$ Since we have a method to compute $m$-th cumulants of $%
x,$ as we described above, we can compute $n$-th moments, by using the M\"{o}%
bius inversion. Actually, we have that

\strut

\begin{center}
$\varphi (x^{n})=\underset{\pi \in NC(n)}{\sum }k_{\pi }\left( \underset{%
n-timse}{\underbrace{x,.....,x}}\right) =\underset{\pi \in NC(n)}{\sum }\,%
\underset{(v_{1},...,v_{k})\in \pi }{\prod }k_{k}\left( \underset{n-times}{%
\underbrace{x,.....,x}}\right) .$
\end{center}

\strut

So,

\strut

\begin{center}
$\varphi (x^{n})=\underset{\pi \in NC(n)}{\sum }\,\underset{%
(v_{1},...,v_{k})\in \pi }{\prod }\left( \underset{\theta \in NC(k)}{\sum }\,%
\underset{(j_{1},...,j_{l})\in \theta }{\prod }\varphi \left(
E(x^{l})\right) \mu (\theta ,1_{k})\right) .$
\end{center}

\strut

\begin{lemma}
Let $B$ be a unital algebra and $A,$ an algebra over $B.$ Suppose that a
NCPSpace $(A,\varphi )$ and an amalgamated NCPSpace over $B,$ $(A,E)$ are
compatible. If $x\in A$ is an operator, then we can compute $n$-th cumulants
($n\in \mathbb{N}$) of $x$ by

\strut

$k_{n}\left( \underset{n-times}{\underbrace{x,.....,x}}\right) =\underset{%
\pi \in NC(n)}{\sum }\,\underset{(v_{1},...,v_{k})\in \pi }{\prod }\varphi
\left( E(x^{k})\right) \mu (\pi ,1_{n}).$

\strut

Also, we can compute $n$-th moments of $x$ by

\strut

$\varphi (x^{n})=\underset{\pi \in NC(n)}{\sum }\,\underset{%
(v_{1},...,v_{k})\in \pi }{\prod }\left( \underset{\theta \in NC(k)}{\sum }\,%
\underset{(j_{1},...,j_{l})\in \theta }{\prod }\varphi \left(
E(x^{l})\right) \mu (\theta ,1_{k})\right) .$

$\square $
\end{lemma}

\strut

By the previous lemma, we can find a scalar-valued R-transform of $x$ and a
scalar-valued moment series of $x,$ under the compatibility ;

\strut

\begin{proposition}
Let $B$ be a unital algebra and $A,$ an algebra over $B.$ Suppose that a
NCPSpace $(A,\varphi )$ and an amalgamated NCPSpace over $B,$ $(A,E)$ are
compatible. If $x\in A$ is an operator, then its scalar-valued R-transform
and moment series as follows ;

\strut

$r_{x}(z)=\sum_{n=1}^{\infty }\left( \underset{\pi \in NC(n)}{\sum }\,%
\underset{(v_{1},...,v_{k})\in \pi }{\prod }\varphi \left(
coef_{k}(M_{x})\right) \mu (\pi ,1_{n})\right) z^{n}$

and

\strut

$m_{x}(z)=\sum_{n=1}^{\infty }\left( \underset{\pi \in NC(n)}{\sum }\,%
\underset{(v_{1},...,v_{k})\in \pi }{\prod }\left( \underset{\theta \in NC(k)%
}{\sum }\,\underset{(j_{1},...,j_{l})\in \theta }{\prod }\varphi \left(
coef_{l}(M_{x})\right) \mu (\theta ,1_{k})\right) \right) z^{n}.$

\strut

In particular, a scalar-valued moment series of $x$ has the following
relation ;

\strut

$m_{x}(z)=r_{x}\left( z\cdot m_{x}(z)\right) .$
\end{proposition}

\strut

\begin{proof}
By the previous lemma, we have, for each $n\in \mathbb{N},$

\strut

$coef_{n}$\strut $\left( r_{x}(z)\right) =k_{n}\left( \underset{n-times}{%
\underbrace{x,.....,x}}\right) $

and

$coef_{n}\left( m_{x}(z)\right) =\varphi (x^{n}).$

\strut

In general, by Nica (See [1]), if $x_{1},...,x_{s}\in (A,\varphi )$ are
scalar-valued random variables ($s\in \mathbb{N}$), then

\strut

$m_{x_{1},...,x_{s}}(z_{1},...,z_{s})=r_{x_{1},...,x_{s}}\left( z_{1}\left(
m_{x_{1},...,x_{s}}(z_{1},...,z_{s})\right) ,...,z_{s}\left(
m_{x_{1},...,x_{s}}(z_{1},...,z_{s})\right) \right) ,$

\strut \strut

since $m_{x_{1},...,x_{s}}=r_{x_{1},...,x_{s}}\,\frame{*}\,\,\,Zeta.$
\end{proof}

\strut \strut

\subsection{Compatible Subalgebras, in $(A,\protect\varphi ),$ relative to $%
(A,E)$}

\strut \strut

\strut

Until now, we observed the compatibility of a NCPSpace $(A,\varphi )$ and an
amalgamated NCPSpace $(A,E).$ It will be interesting to find an ''compatible
part'' in $(A,\varphi ),$ when we have arbitrary $(A,\varphi )$ and $(A,E).$
In this section, we will not assume that $(A,\varphi )$ and $(A,E)$ are
compatible. They are just given NCPSpace and an amalgmated NCPSpace over an
algebra. We want to observe a ''compatible part'' of $(A,\varphi ),$ with
respect to $(A,E).$

\strut

\begin{definition}
Let $B$ be a unital algebra and $A,$ an algebra over $B$ and let $(A,\varphi
)$ be a NCPSpace and $(A,E),$ a NCPSpace over $B.$ A subalgebra $%
A^{o}\subset A$ is called a compatible Subalgebra of $(A,\varphi ),$
relative to $(A,E)$ if

\strut

$\varphi (x)=\varphi \left( E(x)\right) ,$ for all $x\in A^{o}.$
\end{definition}

\strut \strut

\begin{remark}
Let $B$ be a unital algebra and $A,$ an algebra over $B$ and let $(A,\varphi
)$ and $(A,E)$ be given NCPSpace and NCPSpace over $B,$ respectively. Then
there always exists a compatible subalgebra relative to $(A,E)$ called the
scalar-valued compatible subalgebra which is isomorphic to $\mathbb{C}$,

\strut

$\mathbb{C}\cdot 1_{B}=\{\alpha \cdot 1_{B}:\alpha \in \mathbb{C}\}\subset
B\subset A,$

\strut

denoted by $S.$ Clearly, for any $\alpha \cdot 1_{B}=\alpha \cdot 1_{A}\in S$
satisfies

\strut

$\varphi \left( \alpha \cdot 1_{B}\right) =\alpha =\varphi \left( E(\alpha
\cdot 1_{B})\right) .$

\strut

Thus, there always exists compatible subalgebra relative to $(A,E),$ in $%
(A,\varphi ).$
\end{remark}

\strut \strut

\strut The next proposition also shows that there always exists a compatible
subalgebra of a NCPSpace $(A,\varphi )$ relative to an amalgamated NCPSpace $%
(A,E).$

\strut

\begin{proposition}
Let $B$ be a unital algebra and $A,$ an algebra over $B.$ Let $(A,\varphi )$
be a NCPSpace and let $(A,E)$ be a NCPSpace over $B,$ with its $B$%
-functional $E:A\rightarrow B.$ Then $B\subset A$ is a compatible subalgebra
of $A$ relative to $(A,E).$
\end{proposition}

\strut

\begin{proof}
Clearly, we can get that

\strut

$\varphi (b)=\varphi \left( E(b)\right) ,$ for all $b\in B,$

\strut

since $E(b)=b,$ for all $b\in B.$
\end{proof}

\strut \strut

\begin{definition}
\strut Let $B$ be a unital algebra and let $(A,E)$ be a NCPSpace over $B.$
Let $x\in (A,E)$ be a $B$-valued random variable. It is said that a $B$%
-valued random variable $x$ has scalar-valued property if there exists $%
\alpha _{n}\in \mathbb{C}$ such that

\strut

$E(x^{n})=\alpha _{n}\cdot 1_{B}\in B,$

\strut

for all $n\in \mathbb{N}.$
\end{definition}

\strut

\strut Notice that we only use \textbf{trivial} moments of $x,$ not general
moments of $x,$ in the above definition.

\strut \strut

\begin{theorem}
\strut Let $B$ be a unital algebra and $A,$ an algebra over $B.$ Let $%
(A,\varphi )$ be a NCPSpace and $(A,E),$ a NCPSpace over $B,$ with its $B$%
-functional $E:A\rightarrow B.$ Suppose that $x_{0}\in A$ is an operator
which is $B$-central, as a $B$-valued random variable (i.e $x_{0}b=bx_{0},$
for all $b\in B$) and assume that

\strut

$\varphi (x_{0}^{n})=\varphi \left( E(x_{0}^{n})\right) ,$ for all $n\in 
\mathbb{N}.$

\strut

If $x_{0}$ has scalar-valued property determined by a sequence $\left(
\alpha _{n}\right) _{n=1}^{\infty }\subset \mathbb{C}$, then a free product
of an algebra generated by $\{x_{0}\}$ and $B,$ denoted by $A^{o}=A\lg
\left( \{x_{0}\}\right) *B\subset A$, is a compatible subalgebra relative to 
$(A,E). $
\end{theorem}

\strut

\begin{proof}
It suffices to show that

\strut

$\varphi \left(
b_{1}x_{0}^{\,\,\,k_{1}}b_{2}x_{0}^{\,\,\,k_{2}}...b_{n}x_{0}^{\,\,\,k_{n}}%
\right) =\varphi \left(
E(b_{1}x_{0}^{k_{1}}b_{2}x_{0}^{k_{2}}...b_{n}x_{0}^{k_{n}})\right) ,$

\strut

for all $n\in \mathbb{N},$ where $k_{1},...,k_{n}\in \mathbb{N}.$ Remark
that here $\varphi $ can be regarded as $\varphi \mid _{A^{o}},$ where $%
A^{o}=A\lg \left( \{x_{0}\}\right) *B\subset A$. The above equation can be
rewritten by

\strut

$\varphi \left( (b_{1}b_{2}...b_{n})x_{0}^{k_{1}+...+k_{n}}\right) =\varphi
\left( E(b_{1}...b_{n}\cdot x_{0}^{k_{1}+...+k_{n}})\right) $

$\Longleftrightarrow $

$\varphi \left( bx_{0}^{N}\right) =\varphi \left( E(bx_{0}^{N})\right) ,$

\strut

where $b=b_{1}...b_{n}\in B$ and $N=\sum_{j=1}^{n}k_{j},$ since $x_{0}$ is $%
B $-central. Fix $n\in \mathbb{N}.$ By the above equation, it is enough to
show that

\strut

$\varphi \left( bx_{0}^{N}\right) =\varphi \left( bE(x_{0}^{N})\right) ,$
for $b\in B$ arbitrary.

\strut

By scalar-valued property determined by $\left( \alpha _{n}\right)
_{n=1}^{\infty }$ and by the assumption that $\alpha _{n}\cdot
1_{B}=E(x_{0}^{n}),$ for all $n\in \mathbb{N},$

\strut

$\alpha _{n}=\varphi \left( E(x_{0}^{n})\right) ,$ for all $n\in \mathbb{N}.$

\strut And

$\varphi \left( bE(x_{0}^{N})\right) =\alpha _{N}\cdot \varphi (b),$ for all 
$N\in \mathbb{N}.$

\strut

Now, consider that, if we put $y=x_{0}^{N}\in A^{o},$ then

\strut

$\varphi \left( bx_{0}^{N}\right) =\varphi (by)=\underset{\pi \in NC(2)}{%
\sum }k_{\pi }(b,y)$

$\ \ \ \ \ \ \ \ \ \ \ \ \ =k_{2}(b,y)+k_{1}(b)k_{1}(y)=k_{2}(b,y)+\varphi
(b)\varphi (y)$

$\ \ \ \ \ \ \ \ \ \ \ \ \ =0+\varphi (b)\varphi (y)$

\strut

since $b$ and $y$ are free in $(A^{o},\varphi )\equiv \left( A\lg
(\{x_{0}\})*B,\text{ }\varphi \mid _{A\lg (\{x_{0}\})*B}\right) $

\strut

$\ \ \ \ \ \ \ \ \ \ \ \ \ \ =\varphi (b)\varphi \left( x_{0}^{N}\right)
=\varphi (b)\varphi \left( E(x_{0}^{N})\right) $

$\ \ \ \ \ \ \ \ \ \ \ \ \ \ =\varphi (b)\cdot \alpha _{N},$

hence

$\varphi \left( bx_{0}^{N}\right) =\alpha _{N}\cdot \varphi (b),$ for all $%
N\in \mathbb{N}.$

\strut

Therefore, for any $n\in \mathbb{N},$

$\varphi \left( bx_{0}^{n}\right) =\varphi \left( bE(x_{0}^{n})\right) .$
\end{proof}

\strut \strut The above theorem shows how to construct a compatible
subalgebra of $A,$ relative to $(A,E).$ But this is not a general method to
get a compatible subalgebra of $(A,\varphi ),$ relative to $(A,E).$ And this
construction is very artificial. But this provides us one way to construct a
compatible subalgebra containing $B$ in $(A,\varphi ),$ relative to $(A,E).$
We can see the following special case ;

\strut

\begin{example}
Let $B,\,A,\,(A,\varphi )$ and $(A,E)$ be given as before. Let $x_{0}\in A$
be a $B$-semicircular element in $(A,E),$ in the sense that only
nonvanishing $B$-cumulant is the second one, such that

\strut

$K_{2}^{t}(x_{0},x_{0})=t\cdot 1_{B},$ for some $t\in \mathbb{C}.$

and

$x_{0}b=bx_{0},$ \ for all $b\in B.$

\strut

Then $x_{0}$ satisfies the scalar-valued property. i.e, for any $n\in 
\mathbb{N},$

\strut

$E(x_{0}^{n})=0_{B},$ whenever $n$ is odd

and

$E(x_{0}^{2n})=t_{2n}\cdot 1_{B},$ where $t_{2n}\in \mathbb{C}.$

\strut

Indeed, for any $n\in \mathbb{N},$

\strut

$%
\begin{array}{ll}
E\left( x_{0}^{2n}\right) & =\underset{\pi \in NC^{(even)}(2n)}{\sum }\,%
\underset{(v_{1},v_{2})\in \pi }{\prod }K_{2}^{t}(x_{0},x_{0}) \\ 
& =\underset{\pi \in NC^{(even)}(2n)}{\sum }\,\underset{(v_{1},v_{2})\in \pi 
}{\prod }\left( t\cdot 1_{B}\right) \in \mathbb{C}\cdot 1_{B}.%
\end{array}
$

\strut

Remember that a $B$-semicircular element is $B$-even. Hence we only need to
consider $NC^{(even)}(2n).$ Also, recall that, moments of $B$-central
elements are multiplicative, like a scalar-valued case. Hence the second
equality holds true, in this case. (Note that if $x_{0}$ is Not $B$-central,
then we have that

\strut

$E\left( x_{0}^{2n}\right) =\underset{\pi \in NC^{(even)}(2n)}{\sum }\,%
\underset{(v_{1},v_{2})\in \pi (o)}{\prod }\widehat{C}(\pi \mid
_{(v_{1},v_{2})}(x_{0}\otimes \widehat{C}(\pi \mid _{W})(x_{0}\otimes
x_{0})x_{0})),$

\strut

where $W$ is an inner block of $(v_{1},v_{2}),$ if exists.)

\strut

Now, fix a linear functional $\varphi :A\rightarrow \mathbb{C}$. If $x_{0}$
is centered for $\varphi $ and $E$ (i.e $\varphi (x_{0})=0$ and $%
E(x_{0})=0_{B}$), then $\left( A\lg (\{x_{0}\})*B,\,\varphi \right) \subset
(A,\varphi )$ is a compatible subalgebra, relative to $(A,E).$ Since $x_{0}$
is centered for $\varphi $ and $E,$ we have that $\varphi (x_{0})=\varphi
\left( E(x_{0})\right) .$ Moreover, by the previous observation, $\varphi
(x_{0}^{2n})=\varphi \left( E(x_{0}^{2n})\right) ,$ for $n\in \mathbb{N}.$
\end{example}

\strut

\strut \strut

\subsection{Construction of a compatible NCPSpace $(A,\protect\varphi %
^{\prime })$ related to the given $(A,\protect\varphi )$ and $(A,E)$}

\strut

\strut

In this section, we will consider how to construct a NCPSpace, $(A,\varphi
^{\prime }),$ when a NCPSpace $(A,\varphi )$ and an amalgamated NCPSpace $%
(A,E),$ over a unital algebra $B\subset A$ are given. The following
proposition makes us possible to construct a new NCPSpace $(A,\varphi
^{\prime }),$ which is compatible with $(A,E),$ for the given $(A,\varphi ).$

\strut \strut

\begin{proposition}
\strut Let $B$ be a unital algebra and $A,$ an algebra over $B.$ Suppose
that we have a $B$-functional $E:A\rightarrow B$ and a linear functional $%
\psi :B\rightarrow \mathbb{C}.$ Then a NCPSpace $(A,\varphi )$ and a
NCPSpace over $B,$ $(A,E),$ are compatible, where $\varphi =\psi \circ E.$
\end{proposition}

\strut

\begin{proof}
Put $\varphi =\psi \circ E:A\rightarrow \mathbb{C}$. Then it is trivially a
linear functional and hence $(A,\varphi )$ \strut is a NCPSpace. By
definition of $\varphi ,$ we have that

\strut

$\varphi (a)=\psi \left( E(a)\right) =\psi \left( E^{2}(a)\right) =\varphi
\left( E(a)\right) ,$

\strut

for all $a\in A.$ Therefore, $(A,\varphi )$ and $(A,E)$ are compatible via $%
\psi :B\rightarrow \mathbb{C}.$
\end{proof}

\strut \strut

\strut The following theorem shows us that if we have arbitrary NCPSpace $%
(A,\varphi )$ and amalgamated NCPSpace $(A,E)$, where $B\subset A$ and $%
E:A\rightarrow B$ is a $B$-functional, then we can construct a NCPSpace, $%
(A,\varphi ^{\prime }),$ which is compatible with the given amalgamated
NCPSpace, $(A,E)$.

\strut

\begin{theorem}
Let $B$ be a unital algebra and $A,$ an algebra over $B.$ Suppose that we
have a NCPSpace $(A,\varphi )$ and a NCPSpace over $B,$ $(A,E),$ with its $B$%
-functional $E:A\rightarrow B.$ If we define $\varphi ^{\prime }=\varphi
\circ E:A\rightarrow \mathbb{C},$ then a new NCPSpace $(A,\varphi ^{\prime })
$ and $(A,E)$ are compatible.
\end{theorem}

\strut

\begin{proof}
Let $\psi =\varphi \mid _{B}:B\rightarrow \mathbb{C}.$ Then $\psi $ is a
linear functional from $B$ into $\mathbb{C}.$ Define a new linear functional,

\strut

$\varphi ^{\prime }=\psi \circ E:A\rightarrow \mathbb{C}.$

\strut

Then this new linear functional satisfies that

\strut

$\varphi ^{\prime }(a)=\psi \circ E(a)=\varphi \circ E(a)=\varphi \left(
E(a)\right) ,$

\strut

for all $a\in A.$ Therefore, we can get that $(A,\varphi ^{\prime })$ and $%
(A,E)$ are compatible.
\end{proof}

\strut

\strut \strut \strut

\strut

\subsection{Amalgamated Semicircularity and Semicircularity}

\strut \strut

\strut

In this section, we will consider the semicircularity and amalgamated
semicircularity under the compatibility.

\strut

\begin{definition}
(1) Let $(A,\varphi )$ be a NCPSpace. A scalar-valued random variable $a\in
(A,\varphi )$ is semicircular, if only the second (scalar-valued)\ cumulant
is nonvanishing. i.e

\strut

$k_{n}\left( \underset{n-times}{\underbrace{a,.....,a}}\right) =\left\{ 
\begin{array}{lll}
k_{2}(a,a) &  & \text{if }n=2 \\ 
&  &  \\ 
0 &  & \text{otherwise.}%
\end{array}
\right. $

\strut \strut

We say that a family of random variables ($s\in \mathbb{N}$), $%
\{a_{1},...,a_{s}\},$ is a semicircular system if each $a_{j}$ ($j=1,...,s$)
is semicircular and $\{a_{1}\},...,\{a_{s}\}$ are free in $(A,\varphi ).$

\strut

(2) Let $B$ be a unital algebra and $(A,E),$ a NCPSpace over $B,$ with its $%
B $-functional $E:A\rightarrow B.$ A operator-valued random variable $a\in
(A,E)$ is $B$-semicircular (or amalgamated semicircular over $B$) if only
second (operator-valued) cumulant is nonvanishing. i.e

\strut

$%
\begin{array}{ll}
K_{n}\left( \underset{n-times}{\underbrace{a,.....,a}}\right) & 
=C^{(n)}\left( a\otimes b_{2}a\otimes ...\otimes b_{n}a\right) \\ 
&  \\ 
& =\left\{ 
\begin{array}{llll}
C^{(2)}\left( a\otimes b_{2}a\right) =K_{2}(a,a) &  & \text{if }n=2 &  \\ 
&  &  &  \\ 
0_{B} & _{{}} & \text{otherwise}, & 
\end{array}
\right.%
\end{array}
$

\strut

where $b_{2},...,b_{n}\in B$ are arbitrary, for each $n\in \mathbb{N}$ and $%
\widehat{C}=(C^{(n)})_{n=1}^{\infty }\in I\left( A,B\right) $ is the
cumulant multiplicative bimodule map induced by $E,$ in the sense of
Speicher. We say that a family of $B$-valued random variables ($s\in \mathbb{%
N}$), $\{a_{1},...,a_{s}\},$ is a $B$-semicircular family if each $a_{j}$ ($%
j=1,...,s$) is $B$-semicircular and $\{a_{1}\},...,\{a_{s}\}$ are free over $%
B,$ in $(A,E).$
\end{definition}

\strut

\begin{remark}
The above definitions of semircularity and amalgamated semicircularity are
purely combinatorial and purely algebraic. Originally, by Voiculescu, they
are defined on a $C^{*}$- (or $W^{*}$-) algebra framework. Suppose that we
have a unital $C^{*}$-algebra $B$ and $C^{*}$-algebra $A$ containing $B,$ as
its $C^{*}$-subalgebra. Let $\varphi :A\rightarrow \mathbb{C}$ be a state
and $E:A\rightarrow B$, a suitable conditional expectation. ($\varphi
(x^{*})=\overline{\varphi (x)}\in \mathbb{C}$ and $E(x^{*})=E(x)^{*}\in B,$
for all $x\in A$) Then we have a $C^{*}$-probability space $(A,\varphi )$
and an amalgamated $C^{*}$-probability space over $B,$ $(A,E).$ We say that
a scalar-valued random variable $a\in (A,\varphi )$ is semicircular if $a$
is self-adjoint and the only nonvanishing cumulant of $a$ is the second one.
Also, we say that an operator-valued random variable $a\in (A,E)$ is $B$%
-semicircular (or amalgamated semicircular) if $a$ is self-adjoint and the
only nonvanishing operaotr-valued cumulant of $a$ is the second one. So,
originally, we are very much needed the $*$-structure in both scalar-valued
ane operator-valued cases. But, in our definition, we drop this. Of course,
we can consider the $*$-structure on our algebras $B$ and $A$ and we can add
self-adjointness in our definitions. But, in this Section, we will only
consider the combinatorial properties without $*$-structure.
\end{remark}

\strut \strut

If a NCPSpace $(A,\varphi )$ and an amalgamated NCPSpace over $B,$ $(A,E),$
are compatible, there are following semicircular-$B$-semicircular relation ;

\strut

\begin{proposition}
\strut Let $B$ be a unital algebra and $A,$ an algebra over $B.$ Suppose we
have a NCPSpace $(A,\varphi )$ and a NCPSpace over $B,$ $(A,E),$ with its $B$%
-functional $E:A\rightarrow B.$ Assume that $(A,\varphi )$ and $(A,E)$ are
compatible. Then

\strut

(1) if a $B$-central $B$-valued random variable $a\in (A,E)$ is $B$%
-semicircular and if $a\in (A,E)$ satisfies the scalar-valued property , then

\strut

$k_{n}\left( \underset{n-times}{\underbrace{a,.....,a}}\right) =0,$ for any
odd number $n$

and

\strut

$k_{n}\left( \underset{2n-times}{\underbrace{a,.....,a}}\right) =\underset{%
\pi \in NC^{(even)}(2n)}{\sum }\,\left( \underset{(v_{1},...,v_{2k})\in \pi }%
{\prod }\left( \underset{\theta \in NC^{(even)}(2k)}{\sum }\,t^{k(2)}\right)
\right) \mu (\pi ,1_{2n}),$

\strut

where $k(2)$ is the number of pair blocks of $\theta \in NC_{2}(2k),$ only
depending on $k\in \mathbb{N}$, for all $n\in \mathbb{N}$. Here, the numer $%
"t"$ is gotten from

\strut

$K_{2}^{t}(a,a)=t\cdot 1_{B}.$

\strut

(2) if a random variable $a\in (A,\varphi )$ is semicircular and if $a$ is
free from $B,$ in $(A,\varphi ),$ then $a$ is $B$-semicircular, as a $B$%
-valued random variable.
\end{proposition}

\strut

\begin{proof}
(1) Let $a\in (A,E)$ be a $B$-semicircular element. Then, by definition,

\strut

$%
\begin{array}{ll}
K_{n}^{t}\left( \underset{n-times}{\underbrace{a,.....,a}}\right) & 
=C^{(n)}\left( \underset{n-times}{\underbrace{a\otimes .....\otimes a}}%
\right) \\ 
&  \\ 
& =\left\{ 
\begin{array}{lll}
C^{(2)}(a,a)=K_{2}(a,a) &  & \text{if }n=2 \\ 
0_{B} &  & \text{otherwise,}%
\end{array}
\right.%
\end{array}
$

\strut

where $\widehat{C}=(C^{(k)})_{k=1}^{\infty }\in I(A,B)$ is the cumulant
multiplicative bimodule map induced by $E.$ Assume that $K_{2}(a,a)=t\cdot
1_{B}.$ Then,

\strut

$E(a^{n})=0_{B},$ if $n$ is odd

and

$E(a^{2n})=\underset{\pi \in NC^{(even)}(2n)}{\sum }\,\underset{%
(v_{1},v_{2})\in \pi }{\prod }K_{2}(a,a),$

\strut

where

$NC^{(even)}(2n)=\{\theta \in NC(2n):\theta $ does not contain odd blocks$%
\}. $

\strut

(Recall that, by [20], $B$-semicircularity implies the $B$-evenness.)

\strut

Notice that the above equality holds true, because $x_{0}$ is $B$-central
(i.e $x_{0}b=bx_{0},$ $\forall b\in B$). By Section 2.1,

\strut

$k_{n}\left( \underset{n-times}{\underbrace{a,.....,a}}\right) =\underset{%
\pi \in NC(n)}{\sum }\,\underset{(v_{1},...,v_{k})\in \pi }{\prod }\varphi
\left( E(a^{k})\right) \mu (\pi ,1_{n})$

\strut

$\ \ \ \ =\left\{ 
\begin{array}{lll}
0 &  & \text{if }n\text{ is odd} \\ 
&  &  \\ 
\underset{\pi \in NC^{(even)}(n)}{\sum }\,\underset{(v_{1},...,v_{2k})\in
\pi }{\prod }\varphi \left( E(a^{2k})\right) \mu (\pi ,1_{n}) &  & \text{if }%
n\text{ is even.}%
\end{array}
\right. $

\strut

Then, for any $n\in \mathbb{N},$

\strut

$k_{2n}\left( \underset{2n-times}{\underbrace{a,......,a}}\right) =\underset{%
\pi \in NC^{(even)}(2n)}{\sum }\,\left( \underset{(v_{1},...,v_{2k})\in \pi }%
{\prod }\varphi \left( E(a^{2k})\right) \right) \mu (\pi ,1_{2n})$

$\ =\underset{\pi \in NC^{(even)}(2n)}{\sum }\,\left( \underset{%
(v_{1},...,v_{2k})\in \pi }{\prod }\varphi \left( \underset{\theta \in
NC^{(even)}(2k)}{\sum }\,\underset{(w_{1},w_{2})\in \theta }{\prod }%
K_{2}(a,a)\right) \right) \mu (\pi ,1_{2n})$

$\ =\underset{\pi \in NC^{(even)}(2n)}{\sum }\,\left( \underset{%
(v_{1},...,v_{2k})\in \pi }{\prod }\varphi \left( \underset{\theta \in
NC^{(even)}(2k)}{\sum }\,\underset{(w_{1},w_{2})\in \theta }{\prod }t\right)
\right) \mu (\pi ,1_{2n})$

$\ =\underset{\pi \in NC^{(even)}(2n)}{\sum }\,\left( \underset{%
(v_{1},...,v_{2k})\in \pi }{\prod }\left( \underset{\theta \in
NC^{(even)}(2k)}{\sum }\,t^{\left| \theta \right| _{2}}\right) \right) \mu
(\pi ,1_{2n}).$

\strut

(2) Suppose that $a\in (A,\varphi )$ is a semicircular. Since $a$ and $B$
are free in $(A,\varphi ),$ by Section 2.2, we have that

\strut

$%
\begin{array}{ll}
K_{n}^{t}\left( \underset{n-times}{\underbrace{a,......,a}}\right) & 
=k_{n}\left( \underset{n-times}{\underbrace{a,......,a}}\right) \cdot 1_{B}
\\ 
&  \\ 
& =\left\{ 
\begin{array}{lll}
k_{2}(a,a)\cdot 1_{B} &  & \text{if }n=2 \\ 
&  &  \\ 
0_{B} &  & \text{otherwise.}%
\end{array}
\right.%
\end{array}
$

\strut \strut

Therefore, $a$ is $B$-semicircular.
\end{proof}

\strut \strut

\begin{corollary}
Let $B$ be a unital algebra and $A,$ an algebra over $B.$ Suppose that $%
(A,\varphi )$ is a NCPSpace and $(A,E)$ is a NCPSpace over $B,$ with its $B$%
-functional $E:A\rightarrow B.$ Assume that $(A,\varphi )$ and $(A,E)$ are
compatible. Fix $s\in \mathbb{N}.$ If $\{x_{1},...,x_{s}\}$ is a
semicircular system and if $\{x_{1},...,x_{s}\}$ is free from $B,$ in $%
(A,\varphi ),$ then $\{x_{1},...,x_{s}\}$ is a $B$-semicircular system.
\end{corollary}

\strut

\begin{proof}
By the freeness of $\{x_{1},...,x_{s}\}$ and $B,$ by (2) of the previous
proposition, we have that each $x_{j}$ ($j=1,...,s$) is $B$-semicircular,
again. Also, by the compatibility, we have vanishing mixed operator-valued
cumulant of $x_{1},...,x_{s}.$ Indeed, mixed operator-valued cumulants are ;

\strut

$K_{n}\left( x_{i_{1}},...,x_{i_{n}}\right) =C^{(n)}\left( x_{i_{1}}\otimes
b_{i_{2}}x_{i_{2}}\otimes ...\otimes b_{i_{n}}x_{i_{n}}\right) $

\strut

where $b_{i_{2}},...,b_{i_{n}}\in B$ are arbitrary and where $\widehat{C}%
=(C^{(k)})_{k=1}^{\infty }\in I(A,B)$ is the cumulant multiplicative
bimodule map induced by $E$

\strut

$\ \ \ \ \ \ \ \ \ =\left( \varphi (b_{i_{2}})\cdot \cdot \cdot \varphi
(b_{i_{n}})\right) k_{n}^{t}(x_{i_{1}},...,x_{i_{n}})\cdot 1_{B}$

\strut

by the freeness of $\{x_{1},...,x_{s}\}$ and $B,$ in $(A,\varphi )$

\strut

$\ \ \ \ \ \ \ \ \ =\left( \varphi (b_{i_{2}})\cdot \cdot \cdot \varphi
(b_{i_{n}})\right) (0)\cdot 1_{B}=0_{B},$

\strut

by the freeness of $\{x_{1}\},...,\{x_{s}\},$ in $(A,\varphi ).$ Thus $%
\{x_{1},...,x_{s}\}$ is a $B$-semicircular system.
\end{proof}

\strut \strut

\strut \strut

\subsection{Amalgamated Evenness and Evenness}

\strut \strut

\strut

In this section, we will consider the evenness and amalgamated evenness
under the compatibility.

\strut

\begin{definition}
Let $B$ be a unital algebra and $A,$ an algebra over $B.$ Suppose that we
have a NCPSpace $(A,\varphi )$ and an amalgamated NCPSpace over $B,$ $(A,E),$
with its $B$-functional $E:A\rightarrow B.$

\strut

(1) A random variable $a\in (A,\varphi )$ is even if all odd cumulants of $a$
vanish. i.e

\strut

$k_{n}\left( \underset{n-times}{\underbrace{a,......,a}}\right) =0,$ for all
odd number $n\in \mathbb{N}.$

\strut

(2) An operator-valued random variable $a\in (A,E)$ is $B$-even (or
amalgamated even) if all odd cumulants of $a$ vanish. i.e

\strut

$K_{n}\left( \underset{n-times}{\underbrace{a,......,a}}\right)
=C^{(n)}\left( a\otimes b_{2}a\otimes ...\otimes b_{n}a\right) =0_{B},$

\strut

for all odd number $n\in \mathbb{N},$ where $b_{2},...,b_{n}\in B$ are
arbitrary and $\widehat{C}=(C^{(k)})_{k=1}^{\infty }\in I(A,B)$ is the
cumulant multiplicative bimodule map induced by $E.$
\end{definition}

\strut

\begin{remark}
Similar to the semicircular-case, evenness is defined originally on $*$%
-structure. So, self-adjointness is needed. However, like the previous
section, we will only consider the pure combinatorial properties.
\end{remark}

\strut

In [20], we observed that ;

\strut

\begin{proposition}
(See [20]) Let $(A,E)$ be a NCPSpace over a unital algebra $B.$ If a $B$%
-valued random variable $a\in (A,E)$ is $B$-even, then

\strut

(1) all odd moments vanish. (The converse is also true by the M\"{o}bius
inversion.)

\strut

(2) we have that, for any $n\in \mathbb{N},$

\strut

$K_{2n}\left( \underset{2n-times}{\underbrace{a,......,a}}\right)
=C^{(2n)}\left( a\otimes b_{2}a\otimes ...\otimes b_{2n}a\right) $

$\ \ \ \ \ \ \ \ \ \ \ \ \ \ \ \ \ \ \ \ \ \ \ =\underset{\pi \in
NC^{(even)}(2n)}{\sum }\widehat{E}(\pi )\left( a\otimes b_{2}a\otimes
...\otimes b_{2n}a\right) \,\mu (\pi ,1_{2n}),$

\strut

where

$NC^{(even)}(2n)=\{\pi \in NC(2n):\pi $ does not contain odd blocks$\}$

\strut

and $b_{2},...,b_{n}\in B$ are arbitrary. Here, an odd (resp. even) block
means a block containing odd (resp. even) number of entries. $\square $
\end{proposition}

\strut

\begin{proposition}
Let $B,$ $(A,\varphi )$ and $(A,E)$ be given as before. Assume that $%
(A,\varphi )$ and $(A,E)$ are compatible.

\strut

(1) Let $a\in (A,E)$ be even. Then $a\in (A,\varphi )$ is even.

\strut

(2) Let $a\in (A,\varphi )$ be even. If $a$ is free from $B,$ in $(A,\varphi
),$ then $a\in (A,E)$ is $B$-even.
\end{proposition}

\strut

\begin{proof}
(1) Let $a\in (A,E)$ be a $B$-even element. Observe that

\strut

$k_{n}\left( \underset{n-times}{\underbrace{a,......,a}}\right) =\underset{%
\pi \in NC(n)}{\sum }\,\left( \underset{(v_{1},...,v_{k})\in \pi }{\prod }%
\varphi \left( E(a^{k})\right) \right) \mu (\pi ,1_{n})=0,$

\strut

whenever $n$ is odd, since there always exists odd block $%
(v_{1},...,v_{l})\in \pi ,$ for each $\pi \in NC(n)$. So,

\strut

$k_{2n}\left( \underset{2n-times}{\underbrace{a,.....,a}}\right) =\underset{%
\pi \in NC^{(even)}(2n)}{\sum }\,\left( \underset{(v_{1},...,v_{2k})\in \pi }%
{\prod }\varphi \left( E(a^{2k})\right) \right) \mu (\pi ,1_{2n})$

$\ \ \ \ \ \ \ \ \ \ \ \ \ \ \ \ \ \ \ \ \ =\underset{\pi \in NC^{(even)}(2n)%
}{\sum }\,\left( \underset{(v_{1},...,v_{2k})\in \pi }{\prod }\varphi \left(
a^{2k}\right) \right) \mu (\pi ,1_{2n}).$

\strut

Therefore, $a$ is (scalar-valued) even, too.

\strut \strut

(2) Now, suppose that $a\in (A,\varphi )$ is even and assume that $a$ is
free from $B,$ in $(A,\varphi ).$ So, we can get that

\strut

$K_{n}\left( \underset{n-times}{\underbrace{a,......,a}}\right)
=C^{(n)}\left( a\otimes b_{2}a\otimes ...\otimes b_{n}a\right) $

$\strut $

where $b_{2},...,b_{n}\in B$ are arbitrary

\strut

$\ \ \ \ \ \ \ \ \ =\left( \varphi (b_{2})\cdot \cdot \cdot \varphi
(b_{n})\right) k_{n}^{t}\left( \underset{n-times}{\underbrace{a,.....,a}}%
\right) \cdot 1_{B}$

\strut

by the freeness of $a$ and $B,$ in $(A,\varphi )$

\strut

$\ \ \ \ \ \ \ \ \ =0_{B},$

\strut

whenever $n\in \mathbb{N}$ is odd, by the evenness of $a\in (A,\varphi ).$
Therefore, $a$ is $B$-even, in $(A,E).$
\end{proof}

\strut

\begin{theorem}
Let $B$ be a unital algebra and $A,$ an algebra over $B.$ Suppose that a
NCPSpace $(A,\varphi )$ and an amalgamated NCPSpace over $B,$ $(A,E)$ are
compatible. Let $a\in A$ be an operator and let a linear functional $\varphi 
$ be nondegenerated. Then

\strut

$a$ is $B$-even if and only if $a$ is (scalar-valued) even.
\end{theorem}

\strut

\begin{proof}
($\Rightarrow $) It is followed from Proposition 2.18, (1).

\strut \strut

($\Leftarrow $) Assume that $a\in A$ is even. Notice that

\strut

[$a$ is $B$-even]$\overset{def}{\Leftrightarrow }$[all odd $B$-cumulants
vanish]$\overset{(*)}{\Leftrightarrow }$[all odd $B$-moments vanish].

\strut

The (*)-relation holds by the M\"{o}bius inversion (See Proposition 2.7,
(1)). So, it suffices to show that all odd $B$-valued moments of $a$ vanish.
Trivially, since

\strut

$\varphi (x)=\varphi \left( E(x)\right) ,$ for all $x\in A,$

\strut

we have that

$\varphi (a^{n})=\varphi \left( E(a^{n})\right) ,$ for all $n\in \mathbb{N}.$

\strut

Now, let $n$ be odd. Since all odd moments $\varphi (a^{n})$ vanish, all odd 
$B$-moments $E(a^{n})$ vanish, by the nondegeneratedness of $\varphi .$
\end{proof}

\strut

\strut

\strut

\strut \textbf{References}

\strut

\strut

{\small [1] \ \ A. Nica, R-transform in Free Probability, IHP course note,
available at www.math.uwaterloo.ca/\symbol{126}anica.}

{\small [2] \ \ A. Nica, R-transforms of Free Joint Distributions and
Non-crossing Partitions, J. of Func. Anal, 135 (1996), 271-296.}

{\small \strut [3] \ \ A. Nica, R-diagonal Pairs Arising as Free
Off-diagonal Compressions, available at www.math.uwaterloo.ca/\symbol{126}%
anica.}

{\small [4] \ \ A. Nica, D. Shlyakhtenko and R. Speicher, R-cyclic Families
of Matrices in Free Probability, J. of Funct Anal, 188 (2002), 227-271.}

{\small [5] \ \ A. Nica and R. Speicher, R-diagonal Pair-A Common Approach
to Haar Unitaries and Circular Elements, (1995), Preprint}

{\small [6] \ \ A. Nica, D. Shlyakhtenko and R. Speicher, R-diagonal
Elements and Freeness with Amalgamation, Canad. J. Math. Vol 53, Num 2,
(2001) 355-381.}

{\small [7] \ \ A. Nica and R.Speicher, A ''Fouries Transform'' for
Multiplicative Functions on Noncrossing Patitions, J. of Algebraic
Combinatorics, 6, (1997) 141-160.}

{\small [8] \ \ B.Krawczyk and R.Speicher, Combinatorics of Free Cumulants,
available at www.mast.queensu.ca/\symbol{126}speicher.}

{\small \strut \strut [9] \ \ D. Shlyakhtenko, Some Applications of Freeness
with Amalgamation, J. Reine Angew. Math, 500 (1998), 191-212.}

{\small [10] D. Shlyakhtenko, A-Valued Semicircular Systems, J. of Funct
Anal, 166 (1999), 1-47.}

{\small [11] D. Voiculescu, Operations on Certain Non-commuting
Operator-Valued Random Variables, Ast\'{e}risque, 232 (1995), 243-275.}

{\small [12] D.Voiculescu, K. Dykemma and A. Nica, Free Random Variables,
CRM Monograph Series Vol 1 (1992).}

{\small [13] I. Cho, Amalgamated Boxed Convolution and Amalgamated
R-transform Theory (preprint).}

{\small [14] I. Cho, Compressed Amalgamated R-transform Theory, preprint.}

{\small [15] I. Cho, Perturbed R-transform Theory, preprint.}

{\small [16] M. Bo\.{z}ejko, M. Leinert and R. Speicher, Convolution and
Limit Theorems for Conditionally Free Random Variables, Preprint}

{\small [17] P.\'{S}niady and R.Speicher, Continous Family of Invariant
Subspaces for R-diagonal Operators, Invent Math, 146, (2001) 329-363.}

{\small \strut [18] R. Speicher, Combinatorics of Free Probability Theory
IHP course note, available at www.mast.queensu.ca/\symbol{126}speicher.}

{\small [19] R. Speicher, Combinatorial Theory of the Free Product with
Amalgamation and Operator-Valued Free Probability Theory, AMS Mem, Vol 132 ,
Num 627 , (1998).}

{\small [20] R. Speicher, A Conceptual Proof of a Basic Result in the
Combinatorial Approach to Freeness, www.mast.queensu.ca/\symbol{126}speicher.%
}

{\small [21]\ R. Speicher, A. Nica, D. Shlyakhtenko, Operator-Valued
Distributions I. Characterizations of Freeness, preprint.}

\label{REF}

\end{document}